\documentclass[12pt]{article}
\usepackage{mathrsfs}
\usepackage{amsfonts}
\usepackage{amssymb,bbm}  
\usepackage{amsmath,amscd,graphicx,amsthm}
\usepackage{subfigure}
\usepackage{subeqnarray}
\usepackage{cases}
\usepackage{color}
\usepackage{paralist}  
\usepackage[colorlinks,urlcolor=black,anchorcolor=black,linkcolor=black,citecolor=blue]{hyperref}
\usepackage{amsmath}
\usepackage[numbers,sort&compress]{natbib}
\usepackage{indentfirst}
\usepackage{slashed}
\makeatletter \@addtoreset{equation}{section}

       \newtheorem{De}{\bf Definition}[section]
       \newtheorem{lem}{\bf Lemma}[section]
       \newtheorem{remark}{\bf Remark}[section]
       \newtheorem{expl}{\bf Example}[section]
       \newtheorem{assp}{\bf Assumption}

       \newtheorem{thm}{\bf Theorem}[section]

\newcommand{\E}{\mathbb{E}}

\def\be{\begin{equation}}     \def\ee{\end{equation}}
\def\bea{\begin{eqnarray}}    \def\eea{\end{eqnarray}}
\def\beaa{\begin{eqnarray*}}  \def\eeaa{\end{eqnarray*}}

\begin{document}
\title{Propagation of chaos and Razumikhin theorem for the nonlinear McKean-Vlasov SFDEs with common noise}

\author{Xing Chen\thanks{School of Mathematics and Statistics,
Northeast Normal University, Changchun, Jilin, 130024, People's Republic
of China.}
 \and Xiaoyue Li\thanks{School of Mathematical Sciences, Tiangong University, Tianjin, 300387, People's Republic
of China. Research of this author was supported by the National Natural Science Foundation of China (No. 12371402, 12326422), the National Key R$\&$D Program of China (2020YFA0714102), the Natural Science Foundation of Jilin Province, China (No. YDZJ202101ZYTS154) and the Tianjin Natural Science Foundation (24JCZDJC00830)}.
\and Chenggui Yuan\thanks{Department of Mathematics, Swansea University, Bay Campus, Swansea, SA1 8EN, UK.}
}

\date{}

\maketitle
\begin{abstract}
As the limit equations of mean-field particle systems perturbed by common environmental noise, the McKean-Vlasov stochastic differential equations with common noise have received a lot of attention. Moreover, past dependence is an unavoidable natural phenomenon for dynamic systems in life sciences, economics, finance, automatic control, and other fields. Combining the two aspects above, this paper delves into a class of nonlinear McKean-Vlasov stochastic functional differential equations (MV-SFDEs) with common noise. The well-posedness of the nonlinear MV-SFDEs with common noise is first demonstrated through the application of the Banach fixed-point theorem. Secondly, the relationship between the MV-SFDEs with common noise and the corresponding functional particle systems is investigated. More precisely, the conditional propagation of chaos with an explicit convergence rate and the stability equivalence are studied. Furthermore, the exponential stability, an important long-time behavior of the nonlinear MV-SFDEs with common noise, is derived. To this end, the It\^o formula involved with state and measure is developed for the MV-SFDEs with common noise. Using this formula, the Razumikhin theorem is proved, providing an easy-to-implement criterion for the exponential stability. Lastly, an example is provided to illustrate the result of the stability.
\end{abstract}

\section{Introduction}\label{sec:intr}
This paper focuses on a class of nonlinear MV-SFDEs with common noise. Our main motivation stems from two lines of recent advances in the study of stochastic systems and applications. One of them is the emerging interest in McKean-Vlasov stochastic differential equations (MV-SDEs) with common noise, and the other is the unavoidable natural phenomenon of the past dependence in dynamic systems.

Originated from statistical physics, MV-SDEs have attracted significant scrutiny since the pioneering efforts of Kac \cite{K1956} and McKean \cite{Mc1966}. These equations manifest as the limit equations of mean-field particle systems as the number of particles tends to infinite \cite{S1991}. Distinct from the classical It\^o stochastic differential equations, the distributions of the solutions to the MV-SDEs satisfy nonlinear partial differential equations \cite{W2018}. When a common environmental noise perturbs mean-field particle systems, the resulting MV-SDEs with common noise emerge as the limit equations \cite{CD2018}. The conditional distributions of the solutions to the equations with respect to common noise solve the nonlinear stochastic partial differential equations (SPDEs). Due to the inherent connection with both the mean-field particle systems perturbed by a common noise and the nonlinear SPDEs,
the MV-SDEs with common noise have found extensive applications across various domains, including finance, physics, mean-field games, biology and engineering, among others \cite{BLM2023,BC2023,CLZ2023,KT,LL2023,LLW2022,MSZ2018,STW,SW2022,VHP2022}. For the theory of the MV-SDEs with common noise, one can refer to \cite{BSW2023,BLM2023,CD2018,HSS2021,KNRS2022,STW,SW2022,W2021} for the wellposedness, \cite{CD2018,CF2016,STW,SW2022,VHP2022} for the conditional propagation of chaos, \cite{CD2018,CG2019,DTM,KX1999,LAZ2023} for the nonlinear SPDEs, \cite{BW2024,DTM,M,W2021} for the long-time behaviors and references therein.

 Along another line, the presence of delays or past dependence is an unavoidable natural occurrence within dynamic systems in many fields, such as finance, physics, mean-field games, biology, engineering, and so on \cite{BCNPR2019,BYY2016,GN2003,KN1986,SM2017}. Thus, in practical dynamic systems, SFDEs offer a more realistic, efficient, and widely applicable framework than SDEs.  In reference to a diverse range of applications, we point out the studies regarding population systems \cite{NNY2021a,NNY2021b}, stock price models \cite{SM2017}, financial market models \cite{S2005}, the feedback control of Langevin systems \cite{RMT2015}, the feedback control of MV-SDEs \cite{WHGY2022}, among others. Furthermore, recent efforts have been made to the MV-SFDEs, e.g., \cite{BSY2022,HRW2019,HL2023,WXZ2023} and references therein.

 Motivated by the aforementioned two aspects, this paper considers the MV-SFDEs with common noise
  \begin{align}\label{1.1}
\mathrm{d} X(t)
     &=f(X_t ,\mathcal{L}^1( X(t)))
     \mathrm{d} t
      +g( X_t ,\mathcal{L}^1( X(t))) \mathrm{d} B_t
      \nonumber
      \\&\quad+g^0( X_t, \mathcal{L}^1( X(t)))\mathrm{d} B_t^0,~~~~t\geq 0,
      \end{align}
 where  $X_t= \{X(t+\theta),~-\tau\leq \theta\leq 0 \}$ is regarded as a $C$-valued stochastic process, called the  segment of $X $, $f:{C}\times \mathcal{P}(\mathbb{R}^d)\rightarrow \mathbb{R}^d$, $g:{C}\times \mathcal{P}(\mathbb{R}^d)\rightarrow \mathbb{R}^{d\times d}$ and $g^0: {C}\times \mathcal{P}(\mathbb{R}^d)\rightarrow \mathbb{R}^{d\times d}$ are Borel measurable functions.
To the best of our knowledge, there are only a few studies on equations \eqref{1.1} so far. For an example, the well-posedness of linear McKean-Vlasov stochastic delay differential equations with common noise, which is a special case of MV-SFDEs with common noise \eqref{1.1}, is investigated in \cite{CLY2025}. As far as we know, nonlinear MV-SFDEs with the common noise \eqref{1.1} remain unstudied. Therefore, the primary objective of this paper is to delineate the foundational theory of nonlinear MV-SFDEs with common noise \eqref{1.1}, encompassing the well-posedness, conditional propagation of chaos, stability equivalence and stability.

Although models utilizing SFDEs offer greater realism and generality, analyzing these systems, especially their stability, proves to be considerably more difficult. The main difficulty arises from the fact that the space of the functional is an infinite one. Moreover, common noise adds a further dimension of randomness and complexity to the theoretical analysis. Technically, the complex interplay among conditional distributions, conditional expectations, and measure entanglement poses inherent challenges in the research presented in this paper. In the face of these challenges, it is essential to thoroughly understand and treat systems \eqref{1.1} with great care.

For the MV-SFDEs with common noise \eqref{1.1}, we consider the case where the drift term is only locally Lipschitz continuous but satisfies a certain assumption, while the diffusion term is Lipschitz continuous. Borrowing the idea in \cite{CD2018,KNRS2022}, we employ a fixed-point argument to investigate the well-posedness of \eqref{1.1}. The main difficulty here is the construction of a contraction mapping on a suitable space, enabling the use of a fixed-point argument. Taking into consideration the influence of functional terms, we adjust the stochastic process space introduced in \cite{KNRS2022}, and meticulously construct contraction mappings by the conditional distributions with careful derivation. 

 Regarding stability, it is natural to carry over the Lyapunov method to MV-SFDEs with common noise by employing Lyapunov functionals rather than functions. However, constructing Lyapunov functionals is far more difficult than functions, which leads to inconvenience in applications.  
Therefore, taking the influence of conditional distribution into consideration, we aim to establish sufficient criteria for stability by utilizing the variation rate of Lyapunov functions $V(\phi,\mu,t)$ on $\mathbb{R}^d\times\mathcal{P}_2(\mathbb{R}^d)$. To achieve this goal, we  first derive the It\^o formula involved with state and measure for the MV-SFDEs with common noise, enabling us to obtain the derivative of $V$ along the solution to equation \eqref{1.1}. We then investigate stability through the expectation of the derivative of $V$ which is denoted by $\mathbb{E}LV(\phi,\mu,t)$. While it may seem reasonable to require that $\mathbb{E}LV(\phi,\mu,t)$ be negative for all initial data and all $t>0$, such a condition leads to overly stringent restrictions on the coefficients of \eqref{1.1} as discussed in \cite{M2008}. Fortunately, a few moments of reflection in the proper direction indicate that the stability can be guaranteed by only requiring $\mathbb{E}LV(\phi,\mu,t)$ to be negative on some specific local trajectories. This idea was first introduced by Razumikhin \cite{R1956, R1960} for ordinary differential delay equations, and was extended by some researchers to functional differential equations and stochastic functional differential equations (see, for example, \cite{HV1993,M1996}).  Adopting this idea, we develop the Razumikhin theorem for the MV-SDEs with common noise for exponential stability.

The principal contributions of this paper are as follows.
\begin{itemize}
\item  The well-posedness of the nonlinear MV-SFDEs with common noise \eqref{1.1} is proved. We first introduce an appropriate stochastic process space and then construct a contraction mapping by using a frozen SFDE. By applying the Banach fixed point theorem, the existence and uniqueness of the solution to equation \eqref{1.1} in a finite time interval is obtained. Finally, the global solution to equation \eqref{1.1} is finally derived by a splicing method.

\item The relationship between the nonlinear MV-SFDEs with common noise and the mean-field functional particle systems with common noise is established. The conditional propagation of chaos, which describes the convergence, is derived. An explicit bound of the convergence error is given. Then using the convergence, the stability equivalence between the nonlinear MV-SFDEs with common noise \eqref{1.1} and the corresponding particle systems is proved.

\item The stability, an important long-time behavior of the nonlinear MV-SFDEs with common noise, is investigated. For this purpose, we give first the It\^o formula involved with state and measure for the MV-SFDEs with common noise. By this It\^o formula, the Razumikhin theorem is established which gives a criterion based on a Lyapunov function of both state and measure for the exponential stability in the $q$th ($q \geq 2$) moment.
\end{itemize}

The subsequent sections of this paper are organized as follows: Section \ref{Sect.1} gives some preliminaries. Section \ref{Sect.2} investigates the well-posedness of the nonlinear MV-SFDEs with common noise \eqref{1.1}. Section \ref{Sect.3} focuses on the conditional propagation of chaos and the stability equivalence between the equations \eqref{1.1} and the corresponding particle systems. Section \ref{Sect.4} establishes the Razumikhin theorem for the exponential stability of the MV-SFDEs with common noise. Section \ref{Sect.4*} provides an example to illustrate the stability result. Section \ref{Sect.5} concludes this paper.
\section{Preliminaries}\label{Sect.1}
Throughout this paper, we use the following notations.
If $A$ is a matrix or vector, its transpose is denoted by $A^{\prime}$, and its operator norm is denoted by $|A|=\sup\{|Ax|:|x|=1\}=\sqrt{\operatorname{tr}(A^{\prime}A)}$ where $\operatorname{tr}(A^{\prime}A)=\operatorname{trace}(A^{\prime}A)$.
Let $\tau>0$ and denote by ${C}:=C([-\tau,0];\mathbb{R}^d)$ the family of continuous functions $\psi:[-\tau,0]\rightarrow\mathbb{R}^d$ with the norm $\|\psi\|_{{C}}=\sup_{-\tau\leq \theta \leq 0}|\psi(\theta)|$. Let $N\in \mathbb{N}_+$. Let $K$ denote the positive constant which may take different values in different places. We use $K_r$ to emphasize the dependence on $r$.

 Let $(\Omega^0, \mathcal{F}^0, \mathbb{P}^0)$ and $(\Omega^1, \mathcal{F}^1,\mathbb{P}^1)$ be two complete probability spaces equipped with complete and right-continuous filtrations $(\mathcal{F}_t^0)_{t \geq 0}$ and $(\mathcal{F}_t^1)_{t \geq 0}$, respectively. Here $(B^0_t)_{t \geq 0}$ and $(B_t)_{t \geq 0}$  are two independent $d$-dimensional Brownian motions defined on  $\left(\Omega^0, \mathcal{F}^0, \mathbb{P}^0\right)$ and $\left(\Omega^1, \mathcal{F}^1, \mathbb{P}^1\right)$, respectively.  
 Define a product space $(\Omega, \mathcal{F}, \mathbb{P})$, where
$\Omega=\Omega^0 \times \Omega^1,$ {$(\mathcal{F},\mathbb{P})$} is the completion of $(\mathcal{F}^0\otimes \mathcal{F}^1,\mathbb{P}^0 \otimes \mathbb{P}^1)$, and $(\mathcal{F}_t)_{t\geq 0}$ is the complete and right-continuous {augmentation} of $(\mathcal{F}^0_t\otimes\mathcal{F}^1_t)_{t\geq 0}$.
Denote the element of $\Omega$  by $\omega=(\omega^0,\omega^1)$  where $\omega^0\in \Omega^0$ and   $\omega^1\in \Omega^1$.  For $p\geq1$, we use $\mathbb{L}_p(\mathbb{R}^d)$ to denote the family of all $\mathbb{R}^d$-valued random variables $\xi$ such that $\mathbb{E}|\xi|^p<\infty$. Let $\mathbb{L}_p(C)$ be the family of ${C}$-valued random variables $\phi$ satisfying
$\mathbb{E}\|\phi\|_{{C}}^p<\infty$. $\mathcal{P}(\mathbb{R}^d)$ denotes the space of
all probability measures over $\mathbb{R}^d$ equipped with the weak topology. For $p\geq1$, let
$$
\mathcal{P}_p(\mathbb{R}^d):=\Big\{\mu \in \mathcal{P}(\mathbb{R}^d): {\int_{\mathbb{R}^d} |x|^p\mu(\mathrm{d}x) }<\infty\Big\}.
$$
 $\mathcal{P}_p(\mathbb{R}^d)$ is a Polish space under the  Wasserstein distance $\mathbb{W}_p$, where for any $\mu, \nu \in \mathcal{P}_p(\mathbb{R}^d)$,
$$
\mathbb{W}_p(\mu, \nu)\!=\!\inf_{\pi \in \Pi(\mu, \nu)}\Big(\int_{\mathbb{R}^d \times \mathbb{R}^d}|x-y|^p \pi(\mathrm{d} x, \mathrm{~d} y)\Big)^{1/p},
$$
 and $\Pi(\mu, \nu)$ is the set of all couplings of $\mu$ and $\nu$.

For a random variable $\xi$ on $(\Omega, \mathcal{F}, \mathbb{P})$, denote its distribution by $\mathcal{L}(\xi)$, and define a mapping $\mathcal{L}^1(\xi): \Omega^0 \ni\omega^0 \mapsto \mathcal{L}(\xi(\omega^0, \cdot))$. It is known from \cite[Vol-II, Lemma 2.4]{CD2018} that the mapping $\mathcal{L}^1(\xi)$ is a $\mathbb{P}^0$-almost surely well defined and is a random variable from $(\Omega^0, \mathcal{F}^0, \mathbb{P}^0)$  to $\mathcal{P}(\mathbb{R}^d)$ which can be seen as a conditional distribution of $\xi$ with respect to $\mathcal{F}^0\otimes \{\varnothing,\Omega^1\}$. We remark that if the initial data $X_0$ of equation \eqref{1.1} is defined on $(\Omega^1, \mathcal{F}^1_0,\mathbb{P}^1)$, then by \cite[Vol-II, Proposition 2.9]{CD2018}, $\mathcal{L}^1(X(t))$ is a version of the conditional law of $X(t)$ given $B^0$.
 It is worth noting that we do not distinguish a random variable $\xi$ on $(\Omega^0, \mathcal{F}^0, \mathbb{P}^0)$ (resp. $(\Omega^1, \mathcal{F}^1, \mathbb{P}^1)$) from its natural extension $\tilde{\xi}(\omega^0,\omega^1)\equiv \xi(\omega^0)$ (resp. $\tilde{\xi}(\omega^0,\omega^1)\equiv \xi(\omega^1)$) on $(\Omega, \mathcal{F}, \mathbb{P})$.
Denote  the copies of probability space $(\Omega^1, \mathcal{F}^1,\mathbb{P}^1)$ by $(\Omega^1_i, \mathcal{F}^1_i,\mathbb{P}^1_i)$ and the expectation under $\mathbb{P}^1_i$ by $\mathbb{E}^1_i$,  $i=1,2$. Let $\Omega_i=\Omega^0 \times {\Omega}^1_i$,  $i=1,2$.   Let $\delta_{0}$ be the Dirac measure at original point $ 0\in\mathbb{R}^d$. Let $
\varsigma $  be a  $\mathcal{F}_0^1$-measurable $C$-valued random variable on $(\Omega^1, \mathcal{F}^1, \mathbb{P}^1)$.

We now introduce some definitions and some useful Lemmas. We first give the definition of the solution to the MV-SFDE with common noise \eqref{1.1}.
\begin{De}
\rm An $\mathbb{R}^d$-valued stochastic process $X(t)$ on $[-\tau,\infty)$ is called a strong solution to the MV-SFDE with common noise \eqref{1.1} with initial data $
\varsigma\in \mathbb{L}_2(C)$, if it has the following properties:
\begin{itemize}
\item[(i)] It is continuous  and satisfies that for any $T> 0$, $\mathbb{E}[\sup_{-\tau\leq t\leq T}|X(t)|^2]<\infty$. Moreover,
    $(X_t)_{t\geq0}$ is   $ \{\mathcal{F}_t\}_{t\geq0} $-adapted;
\item[(ii)]For any $ T>0$,
     $$
\mathbb{E}\int_0^T\big(|f(X_t,\mathcal{L}^1(X(t)))|+|g( X_t, \mathcal{L}^1(X(t)))|^2
     +|g^0( X_t ,\mathcal{L}^1(X(t)))|^2\big) \mathrm{d} t<\infty;
     $$
\item[(iii)]It together with a continuous $
\big\{\mathcal{F}^0_t\otimes \{\varnothing,\Omega^1\}\big\}_{t\geq0} $-adapted $\mathcal{P}_2(\mathbb{R}^d)$-valued  version of the process $(\mathcal{L}^1(X(t)))_{t\geq0}$  satisfies $\mathbb{P}$-a.s.
    \begin{align*}
X(t)&=\varsigma(0)+\int_0^tf(X_s ,\mathcal{L}^1(X(s)))
     \mathrm{d} s+\int_0^tg(X_s ,\mathcal{L}^1(X(s)))\mathrm{d} B_s
     \\&\quad+\int_0^t g^0(X_s ,\mathcal{L}^1(X(s))) \mathrm{d} B_s^0.
\end{align*}
\end{itemize}
\end{De}

Next, we give the definitions about Lion's derivative.
\begin{De}\label{D1}
\rm
 The operator $\varphi: \mathcal{P}_{2}(\mathbb{R}^d) \rightarrow \mathbb{R}$ is called L-differentiable at $\mu \in \mathcal{P}_{2}(\mathbb{R}^d)$ if there exists a random variable $\xi \in \mathbb{L}_2(\mathbb{R}^d)$ such that
 $\mu=\mathcal{L}(\xi)$ and the lifted function
 $\tilde{\varphi}(\xi):=\varphi(\mathcal{L}(\xi))$ is $\mathrm{Fr\mathrm{\acute{e}}chet}$ differentiable at $\xi$.
 \end{De}

The $\mathrm{Fr\mathrm{\acute{e}}chet}$ derivative $\mathcal{D}\tilde{\varphi}(\xi)$, when seen as an element of $\mathbb{L}_2(\mathbb{R}^d)$
 via the Riesz representation theorem, can be represented as $$\mathcal{D}\tilde{\varphi}(\xi)=\partial_\mu\varphi(\mathcal{L}(\xi))(\xi),$$
 where $\partial_\mu\varphi(\mathcal{L}(\xi)): \mathbb{R}^d\ni y\rightarrow \partial_\mu\varphi(\mathcal{L}(\xi))(y)\in \mathbb{R}^d$ is mean square integral with respect to $\mathcal{L}(\xi)$ (see, e.g., \cite{CCD2022}).
\begin{De}
\rm
 Denote by $\tilde{C}^{2,2,1}:=\tilde{C}^{2,2,1}(\mathbb{R}^d\times\mathcal{P}_2(\mathbb{R}^d )\times[-\tau,\infty);\mathbb{R}_+)$ the family of all operators $\phi=\phi(x,\mu,t):
\mathbb{R}^d \times \mathcal{P}_2(\mathbb{R}^d)\times[-\tau,\infty) \rightarrow\mathbb{R}_+$ satisfying that
\begin{itemize}
\item[(i)]$\phi$ is twice continuously differentiable at $x$ and is continuously differentiable at $t$;
\item[(ii)] For any $(x, t)\in\mathbb{R}^d\times [-\tau,\infty)$, $\phi$ is twice L-differentiable at $\mu$ such that the L-derivative of $\phi$ at $\mu$: $\mathbb{R}^d \times \mathcal{P}_2(\mathbb{R}^d) \times[-\tau,\infty)\times \mathbb{R}^d \ni(x,\mu, t,y) \mapsto \partial_\mu\phi(x,\mu,t)(y) \in \mathbb{R}^d$, the partial derivative of $\partial_\mu\phi(x,\mu,t)(y)$ at $x$: $\mathbb{R}^d\times\mathcal{P}_2(\mathbb{R}^d) \times [-\tau,\infty)\times \mathbb{R}^d \ni(x,\mu,t, y) \mapsto \partial_x\partial_\mu \phi(x,\mu,t)(y) \in \mathbb{R}^d$, the derivative of $\partial_\mu\phi(x,\mu,t)(y)$ at $y$: $\mathbb{R}^d \times \mathcal{P}_2(\mathbb{R}^d) \times[-\tau,\infty) \times \mathbb{R}^d \ni(x,\mu, t,y) \mapsto \partial_y\partial_\mu \phi(x,\mu,t)(y) \in \mathbb{R}^d$   and the L-derivative of $\partial_\mu\phi(x,\mu,t)(y)$ at $\mu$: $\mathbb{R}^d \times \mathcal{P}_2(\mathbb{R}^d) \times[-\tau,\infty)\times \mathbb{R}^d \times \mathbb{R}^d $ $\ni(x,\mu, t,y, z) \mapsto$ $\partial_\mu^2 \phi(x,\mu,t)(y, z)$ $ \in \mathbb{R}^{d\times d}$ all have the versions which are locally bounded and continuous at any points $(x,\mu,t,y)$ with 
    $y\in \mathrm{Supp}(\mu)$ and $(x,\mu,t,y,z)$ with $y,z\in \mathrm{Supp}(\mu)$;
\item[(iii)]For any compact set $\mathcal{K} \subset \mathbb{R}^d \times \mathcal{P}_{2}(\mathbb{R}^d)\times[-\tau,\infty)$,
\begin{align*}
&\sup _{(x,\mu,t) \in \mathcal{K}}\Big[\int_{\mathbb{R}^d}|\partial_\mu \phi(x,\mu,t)(y)|^2 { \mu(\mathrm{d} y )}
+\int_{\mathbb{R}^d}|\partial_y \partial_\mu \phi(x,\!\mu,t)(y)|^2 \mu(\mathrm{d} y)
\\&\quad\quad \quad \quad
+\int_{\mathbb{R}^d} \int_{\mathbb{R}^d}|\partial_\mu^2 \phi(x,\mu,t)(y, z)|^2 \mu(\mathrm{d} y) \mu(\mathrm{d} z)
\\&\quad\quad \quad \quad +\int_{\mathbb{R}^d} |\partial_x \partial_\mu \phi(x,\mu,t)(y)|^2\mu(\mathrm{d} y) \Big]<+\infty .
\end{align*}
\end{itemize}
\end{De}
We then cite a lemma which is useful for the conditional propagation of chaos.
\begin{lem}\label{L*}{\rm
\cite[Theorem 1]{FG2015}} Let $q\geq 1$. Assume that $\rho\in \mathcal{P}_p(\mathbb{R}^d)$ for some $p>q$. Then there exists a positive constant $c$ (only depending on $ d, p, q$) such that, for all $N \geq 1$,
$$
\begin{aligned}
& \mathbb{E}\left(\mathbb{W}^q_q\left(\rho_N, \rho\right)\right) \leq c \mathbb{W}^q_p\left( \rho,\delta_0\right)\varepsilon_N,
\end{aligned}
$$
where
\begin{align}\label{L**}
\varepsilon_N= \begin{cases}N^{-1 / 2}+N^{-(p-q) / p} & \text { if } q>d / 2 \text { and } p \neq 2 q, \\
N^{-1 / 2} \log (1+N)+N^{-(p-q) / p} & \text { if } q=d / 2 \text { and } p \neq 2 q, \\
N^{-q / d}+N^{-(p-q) / p} & \text { if } q\in(0, d / 2) \text { and } p \neq d /(d-q) ,\end{cases}
\end{align}
 and {$\rho_N=\frac{1}{N}\sum_{k=1}^N\delta_{\xi^{k}}$} is the empirical measure of an independent identically distribution sequence {$(\xi^k)_{ k\in \mathbb{N}^*}$} of $\rho$-distributed random variables.
\end{lem}

We end this section with a significant Lemma on the properties of the mapping $\mathcal{L}^1$.
\begin{lem}\label{l2}
 Let $(\xi (t))_{t \geq 0}$ denote a $ \{\mathcal{F}_t\}_{t\geq0} $-adapted $\mathbb{R}^d$-valued process which has continuous paths satisfying $\mathbb{E}[\sup _{0 \leq t \leq T}|\xi(t)|^p]<\infty$ for any $T>0$ and some $ p\geq2$.
Then there exists a version of $(\mathcal{L}^1(\xi(t)) )_{t\geq 0}$, which is $\{\mathcal{F}_t^0\otimes \{\varnothing,\Omega^1\}\}_{t\geq0} $-adapted and has continuous paths in $\mathcal{P}_p(\mathbb{R}^d)$.
 
\end{lem}
{\bf Proof.}
The required assertions follow by the similar way as \cite[Vol-II, Lemma 2.5]{CD2018}. To avoid duplications we omit the proof.
\qed
\section{The wellposedness}\label{Sect.2}
This section gives the well-posedness of the nonlinear MV-SFDE with common noise \eqref{1.1}. For this purpose, we first impose some assumptions. 
\begin{assp}
 \label{A1}
 \rm There exist positive constants  $L_1,L_2$ such that for any $\phi,\psi$ $\in{C}$, and $\mu, \nu \in \mathcal{P}_{2}(\mathbb{R}^d)$,
\begin{align*}
\langle \phi(0)-\psi(0), f( \phi, \mu)-f(\psi, \nu)\rangle&\leq L_1(\|\phi-\psi\|_{{C}}^2+\mathbb{W}_2^2(\mu,\nu)),
\\|g(\phi, \mu)-g(\psi, \nu)|^2
+|g^0(\phi, \mu)-g^0(\psi, \nu)|^2&\leq L_2(\|\phi-\psi\|_{{C}}^2+\mathbb{W}_2^2(\mu,\nu)).
\end{align*}
\end{assp}
\begin{assp}
 \label{A2}
 \rm $f$ is a continuous function on $\mathbb{R}^d\times\mathcal{P}_2(\mathbb{R}^d)$ and is bounded on any bounded set in $\mathbb{R}^d\times\mathcal{P}_2(\mathbb{R}^d)$.
\end{assp}

It follows from Assumption \ref{A1} that there exist positive constants $L_3, L_4$ such that for any $\phi$ $\in{C}$ and $\mu\in \mathcal{P}_{2}(\mathbb{R}^d) $,
\begin{align}
\label{A11*}&\langle \phi(0), f( \phi, \mu)\rangle
\leq L_3(1+\|\phi\|_{{C}}^2+\mathbb{W}_2^2(\mu,\delta_0)),
\\&
\label{A11}|g(\phi, \mu)|^2
+|g^0(\phi, \mu)|^2
\leq L_4(1+\|\phi\|_{{C}}^2+\mathbb{W}_2^2(\mu,\delta_0)).
\end{align}
The main result of this section is given as follows.
\begin{thm}\label{T1}
Let {\rm Assumption \ref{A1}} and {\rm \ref{A2}}
hold. For any  $p\geq 2$ and the initial data $\varsigma\in \mathbb{L}_p(C) $, MV-SFDE with common noise \eqref{1.1} has a uniquely globally strong solution $X(t)$ on  $[-\tau,\infty)$ satisfying that for any $T>0$ $$\mathbb{E}\Big[\sup_{t\in[-\tau,T]}|X(t)|^p\Big]<\infty. $$
\end{thm}
Before proving Theorem \ref{T1}, we do some preparation. We begin with some definitions used later. For $p\geq 2$, 
denote by $\mathbb{L}_1^0( \mathcal{P}_p ):=\mathbb{L}_1(\Omega^0;\mathcal{P}_p(\mathbb{R}^d))$ the family of random variables $\nu :\Omega^0\rightarrow\mathcal{P}_p(\mathbb{R}^d) $ with $\mathbb{E}^0\int_{\mathbb{R}^d}|x|^p\nu(dx)$ $<\infty $. Define $d(\mu,\nu)=(\mathbb{E}^0\mathbb{W}_p^p(\mu,\nu))^{1/p} $, $\mu$, $\nu\in \mathbb{L}_1^0( \mathcal{P}_p )$.
 Due to the completeness of $(\mathcal{P}_p(\mathbb{R}^d),\mathbb{W}_p)$, 
   by virtue of \cite[Lemma A.5]{SS} one notices that $(\mathbb{L}_1^0( \mathcal{P}_p ),d)$ is a complete metric space.
  For any $T>0$, denote by $D ([0,T]) $ the family of continuous maps $\mu$ from $[0,T]$ to $\mathbb{L}_1^0( \mathcal{P}_p )$. Define $\|\mu-\nu\|_{D} =\sup_{0\leq s\leq T}d(\mu(s),\nu(s)),$ $ \mu,\nu\in D ([0,T])$.  Then  $(D ([0,T]),\|\cdot\|_{D })$ is also a complete metric space.

Now we introduce an auxiliary frozen SFDE 
corresponding to \eqref{1.1}.
For any given $\mu\in D ([0,T])$,  replacing  $\mathcal{L}^1$ by $\mu$  in 
\eqref{1.1} yields the frozen SFDE
\begin{align}\label{e2}
y^{\mu}(t)&=y^{\mu}(0)+\int_0^tf( y^{\mu}_s ,\mu_s)
     \mathrm{d} s
+\int_0^tg( y^{\mu}_s ,\mu_s)\mathrm{d} B_s
     +\int_0^t g^0( y^{\mu}_s ,\mu_s)\mathrm{d} B_s^0.
\end{align} 
We now analyze the well-posedness of the frozen SFDE \eqref{e2}.
\begin{lem}\label{l1}
Let {\rm Assumption \ref{A1} }and {\rm \ref{A2}} hold.
For any  $ p\geq 2$ 
and $\mu\in D ([0,T])$,  the equation \eqref{e2} with  the initial data $\varsigma\in \mathbb{L}_p(C)  $ has a uniquely globally strong solution $y^{\mu}({\cdot})$ on $[0,T]$ satisfying
\begin{align}\label{L31}
\mathbb{E}\Big[\sup_{t\in[-\tau,T]}|y^{\mu}(t)|^p\Big] <\infty.
 \end{align}
\end{lem}
{\bf Proof.} By the similar arguments to \cite[Theorem 2.8]{M2008}, one obtains that the frozen SFDE \eqref{e2} has a uniquely globally strong solution $y^{\mu}({\cdot})$ on $[0,T]$ under Assumptions \ref{A1} and \ref{A2}. We now turn to prove \eqref{L31}.
 By the It\^o formula, \eqref{A11*} and \eqref{A11}, one can derive that for any $t\in[0,T]$,
 \begin{align*}
\mathrm{d}[1+|y^{\mu}(t)|^2]^{\frac{p}{2}}
 &=\Big(p[1+|y^{\mu}(t)|^2]^{\frac{p-2}{2}}\langle y^{\mu}(t),f( y^{\mu}_t,\mu_t)\rangle
 \\&\quad+\frac{p(p-1)}{2}[1+|y^{\mu}(t)|^2]^{\frac{p-2}{2}}(|g( y^{\mu}_t,\mu_t)|^2+|g^0( y^{\mu}_t,\mu_t)|^2)\Big) \mathrm{d}t
 \\&
\quad+p[1\!+\!|y^{\mu}(t)|^2]^{\frac{p-2}{2}}\langle y^{\mu}(t),g( y^{\mu}_t,\mu_t)\rangle \mathrm{d}B_t
  \\&\quad+\!p[1\!+\!|y^{\mu}(t)|^2]^{\frac{p-2}{2}}\langle y^{\mu}(t),g^0( y^{\mu}_t,\mu_t)\rangle \mathrm{d}B^0_t
\\
  & \leq\Big(\frac{2L_3 p \!+\!L_4p(p\!-\!1)}{2}[1\!+\!|y^{\mu}(t)|^2]^{\frac{p-2}{2}}(1\!+\! \|y^{\mu}_t\|^2_C\!+\!\mathbb{W}_2^2(\mu_t,\delta_0))\Big) \mathrm{d}t
\\
 &\quad+p[1\!+\!|y^{\mu}(t)|^2]^{\frac{p-2}{2}}\langle y^{\mu}(t),g( y^{\mu}_t,\mu_t)\rangle \mathrm{d}B_t
\\&
 \quad+p[1\!+\!|y^{\mu}(t)|^2]^{\frac{p-2}{2}}\langle y^{\mu}(t),g^0( y^{\mu}_t,\mu_t)\rangle \mathrm{d}B^0_t.
 \end{align*}
Using the Young inequality and the discrete H\"older inequality (see, for instance, \cite[page 52, inequality (2.7) and (2.8)]{MY2006}) shows that
 \begin{align}\label{eq3.5}
 &\mathbb{E}\Big[\sup_{s\in[0,t]}[1+|y^{\mu}(s)|^2]^{\frac{p}{2}}\Big]
 \nonumber\\& \leq2^{\frac{p-2}{2}}(1+\mathbb{E}\|\varsigma\|^p_{C})
 +\big(2L_3+L_4(p-1)\big)\mathbb{E}\int_0^t\mathbb{W}_2^p(\mu_s,\delta_0)
 \mathrm{d}s
 \nonumber  \\
 &\quad
 +{\big(2L_3 +L_4(p-1)\big)\big(p-1\big)}\mathbb{E}\int_0^t[1+\|y^{\mu}_s\|_C^2]^{\frac{p}{2}}\mathrm{d}s
 \nonumber  \\
 &\quad+p\mathbb{E}\Big(\sup_{s\in[0,t]}\int_0^s[1+|y^{\mu}(r)|^2]^{\frac{p-2}{2}}\langle y^{\mu}(r),g( y^{\mu}_r,\mu_r)\rangle \mathrm{d}B_r\Big)
\nonumber \\&\quad
 +p \mathbb{E}\Big(\sup_{s\in[0,t]}\int_0^s[1+|y^{\mu}(r)|^2]^{\frac{p-2}{2}}\langle y^{\mu}(r),g^0( y^{\mu}_r,\mu_r)\rangle \mathrm{d}B^0_r\Big).
 \end{align}
Making use of the Burkholder-Davis-Gundy inequality, the Young inequality and \eqref{A11}, we obtain that
\begin{align}\label{eq3.6}
&p\mathbb{E}\Big(\sup_{s\in[0,t]}\int_0^r[1+|y^{\mu}(r)|^2]^{\frac{p-2}{2}}\langle y^{\mu}(r),g( y^{\mu}_r,\mu_r)\rangle \mathrm{d}B_r\Big)
\nonumber\\&\leq 4\sqrt{2}p\mathbb{E}\Big(\int_0^t[1+|y^{\mu}(s)|^2]^{{p-2}}| y^{\mu}(s)g( y^{\mu}_s,\mu_s)|^2 \mathrm{d}s\Big)^{\frac{1}{2}}
\nonumber\\&\leq 4\sqrt{2}p\mathbb{E}\Big\{\Big(\sup_{s\in[0,t]}[1+|y^{\mu}(s)|^2]^{{\frac{p}{2}}}\Big)\int_0^t[1+|y^{\mu}(s)|^2]^{\frac{{p-2}}{2}}|g( y^{\mu}_s,\mu_s)|^2 \mathrm{d}s\Big\}^{\frac{1}{2}}
\nonumber\\&
\leq\frac{1}{4}\mathbb{E}\Big(\sup_{s\in[0,t]}[1+|y^{\mu}(s)|^2]^{{\frac{p}{2}}}\Big)+64{p}^2\mathbb{E}\int_0^t[1+|y^{\mu}(s)|^2]^{\frac{{p-2}}{2}}|g( y^{\mu}_s,\mu_s)|^2 \mathrm{d}s
\nonumber\\
&
\leq\frac{1}{4}\mathbb{E}\Big(\sup_{s\in[0,t]}[1+|y^{\mu}(s)|^2]^{{\frac{p}{2}}}\Big)+64L_4{p}^2\mathbb{E}\int_0^t[1+|y^{\mu}(s)|^2]^{\frac{{p-2}}{2}}
\nonumber\\&\quad \times (1+\|y^{\mu}_s\|_C^2+\mathbb{W}_2^2(\mu_s,\delta_0))\mathrm{d}s
\nonumber\\
&\leq\frac{1}{4}\mathbb{E}\Big(\sup_{s\in[0,t]}[1+|y^{\mu}(s)|^2]^{{\frac{p}{2}}}\Big)+128L_4p(p-1)\mathbb{E}\int_0^t[1+\|y^{\mu}_s\|_C^2]^{\frac{{p}}{2}}\mathrm{d}s
\nonumber\\&\quad+128L_4p\mathbb{E}\int_0^t\mathbb{W}_2^p(\mu_s,\delta_0)\mathrm{d}s.
\end{align}
Similarly,
\begin{align}\label{eq3.7}
\begin{aligned}
&p\mathbb{E}\Big(\sup_{s\in[0,t]}\int_0^r[1+|y^{\mu}(r)|^2]^{\frac{p-2}{2}}\langle y^{\mu}(r),g^0( y^{\mu}_r,\mu_r)\rangle \mathrm{d}B_r\Big)
\\%
&\leq\frac{1}{4}\mathbb{E}\Big(\sup_{s\in[0,t]}[1+|y^{\mu}(s)|^2]^{{\frac{p}{2}}}\Big)+128L_4p(p-1)\mathbb{E}\int_0^t[1+\|y^{\mu}_s\|_C^2]^{\frac{{p}}{2}}\mathrm{d}s
\\&\quad
+128L_4p\mathbb{E}\int_0^t\mathbb{W}_2^p(\mu_s,\delta_0)\mathrm{d}s.
\end{aligned}\end{align}
Combining \eqref{eq3.5}, \eqref{eq3.6} and \eqref{eq3.7} gives
 \begin{align*}
 \mathbb{E}\Big[\sup_{s\in[0,t]}[1+|y^{\mu}(s)|^2]^{\frac{p}{2}}\Big]&\leq2^{\frac{p-2}{2}}(1+\mathbb{E}\|\varsigma\|^p_{C})+\frac{1}{2}\mathbb{E}\Big(\sup_{s\in[0,t]}[1+|y^{\mu}(s)|^2]^{{\frac{p}{2}}}\Big)
 \end{align*}
 \begin{align*}
 &\quad+\big(2L_3 +L_4(257p-1)\big)\big(p-1\big)\mathbb{E}\int_0^t[1+\|y^{\mu}_s\|_C^2]^{\frac{p}{2}}\mathrm{d}s
\\&\quad+\big(2L_3+L_4(257p-1)\big)\mathbb{E}\int_0^t\mathbb{W}_2^p(\mu_s,\delta_0)
 \mathrm{d}s.
 \end{align*}
 Therefore,
  \begin{align}\label{eq3.8}
 \begin{aligned}
 &\mathbb{E}\Big[\sup_{s\in[0,t]}[1+|y^{\mu}(s)|^2]^{\frac{p}{2}}\Big]
 \\&\leq2^{\frac{p}{2}}(1+\mathbb{E}\|\varsigma\|^p_{C})
 +K\int_0^t\mathbb{E}\Big[[1+\|y^{\mu}_s\|_C^2]^{\frac{p}{2}}\Big]\mathrm{d}s
 +K\mathbb{E}\int_0^t\mathbb{W}_2^p(\mu_s,\delta_0)\mathrm{d}s,
  \end{aligned}
 \end{align}
 where $K=2\big(2L_3 +L_4(257p-1)\big)\big(p-1\big)$.
 Note that
 \begin{align*}
\mathbb{E}\Big[\sup_{s\in[-\tau,t]}[1+|y^{\mu}(s)|^2]^{\frac{p}{2}}\Big]
 &\leq \mathbb{E}[1+\|\varsigma\|_C^2]^{\frac{p}{2}}+\mathbb{E}\Big[\sup_{s\in[0,t]}[1+|y^{\mu}(s)|^2]^{\frac{p}{2}}\Big]
 \\&\leq 2^{\frac{p-2}{2}}(1+\mathbb{E}\|\varsigma\|_C^p)+\mathbb{E}\Big[\sup_{s\in[0,t]}[1+|y^{\mu}(s)|^2]^{\frac{p}{2}}\Big].
 \end{align*}
 Then it follows from \eqref{eq3.8} and the H\"older inequality that
  \begin{align*}
 &\mathbb{E}\Big[\sup_{s\in[-\tau,t]}[1+|y^{\mu}(s)|^2]^{\frac{p}{2}}\Big]
 \\& \leq\frac{3}{2}2^{\frac{p}{2}}(1+\mathbb{E}\|\varsigma\|^p_{C})
 +K\int_0^t\mathbb{E}\Big[\sup_{r\in[0,s]}[1+\|y^{\mu}_r\|_C^2]^{\frac{p}{2}}\Big]\mathrm{d}s
 +K\mathbb{E}\int_0^t\mathbb{W}_2^p(\mu_s,\delta_0)\mathrm{d}s
  \\& \leq\frac{3}{2}2^{\frac{p}{2}}(1+\mathbb{E}\|\varsigma\|^p_{C})
 +K\int_0^t\mathbb{E}\Big[\sup_{r\in[-\tau,s]}[1+|y^{\mu}(r)|^2]^{\frac{p}{2}}\Big]\mathrm{d}s
 +K\mathbb{E}\int_0^t\mathbb{W}_p^p(\mu_s,\delta_0)\mathrm{d}s.
 \end{align*}
 Then utilizing the Gronwall inequality and the fact that $\mu_{\cdot}\in D ([0,T])$, we derive that
  \begin{align*}
 \begin{aligned}
\mathbb{E}\Big[\sup_{s\in[-\tau,t]}[1+|y^{\mu}(s)|^2]^{\frac{p}{2}}\Big]
\leq\Big(\frac{3}{2}2^{\frac{p}{2}}(1+\mathbb{E}\|\varsigma\|^p_{C})+Kt\sup_{0\leq s\leq t}\mathbb{E}[\mathbb{W}_p^p(\mu_s,\delta_0)]\Big)e^{Kt}<\infty.
  \end{aligned}
 \end{align*}
Then we have
\begin{align*}
\mathbb{E}\Big[\sup_{t\in[-\tau,T]}|y^{\mu}(t)|^p\Big]\leq\mathbb{E}\Big[\sup_{t\in[-\tau,T]}[1+|y^{\mu}(t)|^2]^{\frac{p}{2}}\Big]<\infty.
\end{align*}
 Therefore the proof is complete.
\qed

We are now in the position to prove the well-posedness of the MV-SFDE with common noise \eqref{1.1}.

{\bf Proof of Theorem 3.1.} Borrowing the idea from \cite{CD2018}, we employ the Banach fixed point theorem to establish the existence and uniqueness of the solution for the MV-SFDE with common noise \eqref{1.1}. For any $T>0$, define $D_{\varsigma}([0,T])=\{\mu\in D ([0,T]): \mu_0=\mathcal{L}^1(\varsigma(0))\}$. Then it is easy to obtain that $(D_{\varsigma}([0,T]),\|\cdot\|_{D})$ is a complete metric space. For any $\mu \in D_{\varsigma}([0,T])$,  by using Lemma \ref{l2} and Lemma \ref{l1}, one notes that the conditional distribution $\mathcal{L}^1(y^\mu):=\big(\mathcal{L}^1(y^\mu(t))\big)_{t\in[0,T]}$ of the frozen SFDE \eqref{e2} is $\{\mathcal{F}_t^0\otimes \{\varnothing,\Omega^1\}\}_{t\in[0,T]} $-adapted and satisfies that $\mathcal{L}^1(y^\mu)\in D_{\varsigma}([0,T] )$. Define an operator
$\Phi:D_{\varsigma}([0,T] ) \rightarrow D_{\varsigma}([0,T]  )$ by
$\Phi(\mu) = \mathcal{L}^1(y^\mu ),$ which  is well-defined. If the operator $\Phi$ has a unique fixed point $\mu^0$ on the Banach space $D_{\varsigma}([0,T]  )$, then 
 $y^{\mu^0}(t)$ is the unique strong solution of \eqref{1.1} on $[0,T]$.
 It is therefore sufficient to demonstrate the well-posedness of \eqref{1.1} by proving that the operator
$\Phi$ possesses a unique fixed point.
Let $\mu,\nu\in D_{\varsigma}([0,T] )$.
Making use of the It\^o formula, we arrive at
\begin{align*}
&\mathrm{d}|y^{\mu}(t)-y^{\nu}(t)|^p
 \\&=\Big(p|y^{\mu}(t)-y^{\nu}(t)|^{{p-2}}\langle y^{\mu}(t)-y^{\nu}(t),f( y^{\mu}_t,\mu_t)-f( y^{\nu}_t,\nu_t)\rangle \\&\quad+\frac{p(p-1)}{2}|y^{\mu}(t)-y^{\nu}(t)|^{{p-2}}(|g( y^{\mu}_t,\mu_t)-g( y^{\nu}_t,\nu_t)|^2
 \\&\quad +|g^0( y^{\mu}_t,\mu_t)-g^0( y^{\nu}_t,\nu_t)|^2)\Big) \mathrm{d}t
   \\&\quad+p|y^{\mu}(t)-y^{\nu}(t)|^{{p-2}}\langle y^{\mu}(t)-y^{\nu}(t),g( y^{\mu}_t,\mu_t)-g( y^{\nu}_t,\nu_t)\rangle \mathrm{d}B_t
  \\&\quad+p|y^{\mu}(t)-y^{\nu}(t)|^{{p-2}}\langle y^{\mu}(t)-y^{\nu}(t),g^0( y^{\mu}_t,\mu_t)-g^0( y^{\nu}_t,\nu_t)\rangle \mathrm{d}B^0_t.
  \end{align*}
Applying Assumption \ref{A1} and the Young inequality gives that for any $t\in[0,T]$,
\begin{align}\label{eq3.10}
&\mathbb{E}\Big[\sup_{s\in[0,t]}|y^{\mu}(s)\!-\!y^{\nu}(s)|^p\Big]
\nonumber\\&\leq\Big(\frac{2L_1 p\!+\!L_2 p(p\!-\!1)}{2}\Big)\mathbb{E}\int_0^t|y^{\mu}(s)\!-\!y^{\nu}(s)|^{p-2}(\|y^{\mu}_s\!-\!y^{\nu}_s\|_C^{2}\!+\!\mathbb{W}_2^2(\mu_s,\nu_s))\mathrm{d}s
\nonumber\\
&\quad+J_1(t)+J_2(t)
\nonumber\\
& \leq\mathbb{E}\int_0^t\!\!\Big(\big(2L_1\!+\!L_2(p\!-\!1)\big)\big(p\!-\!1\big)\|y^{\mu}_s\!-\!y^{\nu}_s\|_C^{p}\!+\!\big(2L_1\!+\!{L_2(p\!-\!1)}\big)\mathbb{W}_2^p(\mu_s,\nu_s))\!\Big) \mathrm{d}s
\nonumber\\
&\quad+J_1(t)+J_2(t)
\end{align}
where \begin{align*}&J_1(t)=\mathbb{E}\Big[\sup_{s\in[0,t]}\int_0^s\!\!\!\big(p|y^{\mu}(r)\!-\!y^{\nu}(r)|^{{p-2}}\langle y^{\mu}(r)\!-\!y^{\nu}(r),g( y^{\mu}_r,\mu_r)\!-\!g( y^{\nu}_r,\nu_r)\rangle\big) \mathrm{d}B_r\!\Big],\\ &J_2(t)=\mathbb{E}\Big[\sup_{s\in[0,t]}\int_0^s\!\!\!\big(p|y^{\mu}(r)\!-\!y^{\nu}(r)|^{{p-2}}\langle y^{\mu}(r)\!-\!y^{\nu}(r),g^0( y^{\mu}_r,\mu_r)\!-\!g^0( y^{\nu}_r,\nu_r)\rangle\big) \mathrm{d}B^0_r\!\Big].\end{align*}
Using the Burkholder-Davis-Gundy inequality, the Young inequality, the discrete H\"older inequality and Assumption \ref{A1}, we arrive at
\begin{align} \label{eq3.11}%
  J_1(t)
&\leq 4\sqrt{2}p\mathbb{E}\Big[\sup_{s\in[0,t]}|y^{\mu}(s)\!-\!y^{\nu}(s)|^{{p}}
\int_0^t\!\!\!|y^{\mu}(s)\!-\!y^{\nu}(s)|^{p-2}
  |g( y^{\mu}_s,\mu_s)\!-\!g( y^{\nu}_s,\nu_s)|^2 \mathrm{d}s\Big]^\frac{1}{2}
\nonumber\\
&\leq\frac{1}{4}\mathbb{E}\Big[\!\sup_{s\in[0,t]}|y^{\mu}(s)\!-\!y^{\nu}(s)|^p\Big]\!
+\!64p^2\mathbb{E}\!\!\int_0^t\!\!\!|y^{\mu}(s)\!-\!y^{\nu}(s)|^{p-2}
|g( y^{\mu}_s,\mu_s)\!-\!g( y^{\nu}_s,\nu_s)|^2 \mathrm{d}s
\nonumber \\
  &\leq\frac{1}{4}\mathbb{E}\Big[\sup_{s\in[0,t]}|y^{\mu}(s)\!-\!y^{\nu}(s)|^p\Big]
\!
+\!64L_2 p^2 \mathbb{E}\!\!\int_0^t\!\!\!|y^{\mu}(s)\!-\!y^{\nu}(s)|^{p-2}
\nonumber\\&\quad\times(\|y^{\mu}_s\!-\!y^{\nu}_s\|_C^2\!+\!\mathbb{W}_2^2(\mu_s,\nu_s)) \mathrm{d}s
%
\nonumber\\
&
\leq\frac{1}{4}\mathbb{E}\Big[\sup_{s\in[0,t]}|y^{\mu}(s)\!-\!y^{\nu}(s)|^p\Big]+\!128L_2p(p-1) \mathbb{E}\!\!\int_0^t\!\!\!\|y^{\mu}_s\!-\!y^{\nu}_s\|_C^{p}\mathrm{d}s \nonumber\\&\quad
+\!128L_2p\mathbb{E}\int_0^t\mathbb{W}_2^p(\mu_s,\nu_s) \mathrm{d}s.
\end{align}
Similarly, we have
\begin{align}\label{eq3.12}
  J_2(t)&\leq\frac{1}{4}\mathbb{E}\Big[\sup_{s\in[0,t]}|y^{\mu}(s)-y^{\nu}(s)|^p\Big]+128L_2p(p-1) \mathbb{E}\int_0^t\|y^{\mu}_s-y^{\nu}_s\|_C^{p}\mathrm{d}s
   \nonumber\\&\quad+128L_2p\mathbb{E}\int_0^t\mathbb{W}_2^p(\mu_s,\nu_s) \mathrm{d}s.
\end{align}
Inserting \eqref{eq3.11} and \eqref{eq3.12} into \eqref{eq3.10}, one obtains that
\begin{align*}
&\mathbb{E}\Big[\sup_{s\in[0,t]}|y^{\mu}(s)-y^{\nu}(s)|^p\Big] 
\nonumber \\
&\leq\frac{1}{2}\mathbb{E}\Big[\sup_{s\in[0,t]}|y^{\mu}(s)\!-\!y^{\nu}(s)|^p\Big]+\big(2L_1+L_2(257 p-1)\big)\big(p-1\big)\mathbb{E}\int_0^t\|y^{\mu}_s-y^{\nu}_s\|_C^{p}\mathrm{d}s
\nonumber \\&\quad+\big(2L_1+{L_2(257p-1)}\big)\mathbb{E}\int_0^t\mathbb{W}_2^p(\mu_s,\nu_s) \mathrm{d}s.
\end{align*}
Therefore, the Fubini theorem derives that
\begin{align*}
&\mathbb{E}\Big[\sup_{s\in[0,t]}|y^{\mu}(s)-y^{\nu}(s)|^p\Big]
\\
&\leq K_p\int_0^t\mathbb{E}\Big[\sup_{r\in[0,s]}\|y^{\mu}_r-y^{\nu}_r\|_C^{p}\Big]\mathrm{d}s
+K_p\mathbb{E}\int_0^t\mathbb{W}_2^p(\mu_s,\nu_s) \mathrm{d}s
\\&
\leq K_p\int_0^t\mathbb{E}\Big[\sup_{r\in[0,s]}|y^{\mu}(r)-y^{\nu}(r)|^{p}\Big]\mathrm{d}s
+K_p\mathbb{E}\int_0^t\mathbb{W}_2^p(\mu_s,\nu_s) \mathrm{d}s,
\end{align*}
where $K_p=2\big(2L_1+L_2(257 p-1)\big)\big(p-1\big)$.
Using the Gronwall inequality and the H\"older inequality gives that for any $t\in [0, T]$,
\begin{align}\label{eq3.13}
\mathbb{E}\Big[\sup_{s\in[0,t]}|y^{\mu}(s)-y^{\nu}(s)|^p\Big]
&\leq  K_p e^{K_{p}T}\int_0^t\sup_{r\in[0,s]}
\mathbb{E}\mathbb{W}_2^p(\mu_r,\nu_r)\mathrm{d} r
\nonumber\\&\leq K_p e^{K_{p}T}\int_0^t\sup_{r\in[0,s]}
\mathbb{E}\mathbb{W}_p^p(\mu_r,\nu_r)\mathrm{d} r.
\end{align}
One derives from the property of expectation and the fact that for any $t\in[0,T]$, $\mathbb{E}^0\mathbb{W}_p^p(\mu_t,\nu_t)$ is a constant that
\begin{align}\label{eq3.13*}
\mathbb{E}\mathbb{W}_p^p(\mu_t,\nu_t)=\mathbb{E}[\mathbb{E}^0\mathbb{W}_p^p(\mu_t,\nu_t)]=\mathbb{E}^0\mathbb{W}_p^p(\mu_t,\nu_t).
\end{align}
Combining \eqref{eq3.13} and \eqref{eq3.13*} shows that for any $t\in [0, T]$,
\begin{align*}
\mathbb{E}\Big[\sup_{s\in[0,t]}|y^{\mu}(s)-y^{\nu}(s)|^p\Big]
&\leq K_p e^{K_{p}T}\int_0^t\sup_{r\in[0,s]}
\mathbb{E}^0\mathbb{W}_p^p(\mu_r,\nu_r)\mathrm{d} r
\\&\leq K_pT e^{K_{p}T}\sup_{s\in[0,T]}
\mathbb{E}^0\mathbb{W}_p^p(\mu_s,\nu_s)\\&
\leq  {K_T}\|\mu-\nu\|^p_{D},
\end{align*}
where $K_T={K_{p}Te^{K_{p}T}}$. This, together with the metric definition on $D ([0,T])$, implies
\begin{align*}
 \|\Phi(\mu)-\Phi(\nu)\|^p_{D}&=\sup_{t\in[0,T]}\mathbb{E}^0\mathbb{W}_p^p(\mathcal{L}^1(y^{\mu}(t)),\mathcal{L}^1(y^{\nu}(t)))
\\
&\leq\sup_{t\in[0,T]}\mathbb{E}^0[\mathbb{E}^1|y^{\mu}(t)-y^{\nu}(t)|^p]
 \\&=\sup_{t\in[0,T]}\mathbb{E}[|y^{\mu}(t)-y^{\nu}(t)|^p]
 \\&\leq \mathbb{E}\Big[\sup_{t\in[0,T]}|y^{\mu}(t)-y^{\nu}(t)|^p\Big]
 \leq {K_T}\|\mu-\nu\|^p_{D}.
\end{align*}
 Choose a constant  $T_1>0$ small sufficiently such that $K_{T_1}<1$. It follows from the Banach fixed point theorem  (see, for instance, \cite[Theorem 2.2]{A2014}) that  $\Phi$ owns a unique fixed point on $D_{\varsigma}([0,T_1])$. This implies that the equation \eqref{1.1} exists a solution on $[0,T_1]$. It is easy to see that a solution $X(t)$ of the equation \eqref{1.1} must be a solution of the equation \eqref{e2}  with $\mu_{t}=\mathcal{L}^1(X(t))$ on $[0,T_1]$. Since the equation \eqref{e2} has a unique solution under Assumption \ref{A1} and \ref{A2}, the equation \eqref{1.1} with the initial data $\varsigma$ has a unique strong solution $X(t)$ on $[0,T_1]$ satisfying $\mathbb{E}[\sup_{t\in[-\tau,T_1]}|X(t)|^p]< \infty$. Next, consider the equation \eqref{1.1} on $[T_1, 2T_1]$ with initial data $X_{T_1}$. By looking for the unique fixed point of mapping $\Phi$ on $D_{X_{T_1}}([T_1, 2T_1] )$, we obtain the unique strong solution $X(t)$ of \eqref{1.1} on $[T_1 ,2T_1]$ satisfying $\mathbb{E}[\sup_{t\in[T_1-\tau ,2T_1]}|X(t)|^p]<\infty$. The desired assertions follow by repeating this procedure on $[2T_1, 3T_1],~\cdots$. \qed
  \section{Conditional propagation of chaos and the stability equivalence}\label{Sect.3}
This section firstly investigates the conditional propagation of chaos which describes the convergence between the equation \eqref{1.1} and the corresponding particle systems with common noise. Then the stability equivalence between the equation \eqref{1.1} and the corresponding particle systems with common noise is established. Let $(\varsigma^k, (B^k_t)_{t \geq 0})_{k\geq1}$ be a sequence of independent copies of $(\varsigma, (B_t)_{t \geq 0})$. For each integer $k\geq 1$, let $\bar{X}^{k}(t)$ be the solution to the MV-SFDE with common noise
 \begin{align}\label{eq4.2}
\begin{aligned}
\mathrm{d}\bar{X}^{k}(t)&= f(\bar{X}^{k}_t,\mathcal{L}^1(\bar{X}^{k}(t)))\mathrm{d} t
+g(\bar{X}^{k}_t,\mathcal{L}^1(\bar{X}^{k}(t)))\mathrm{d}B^k_t
\\&\quad
+g^0(\bar{X}^{k}_t,\mathcal{L}^1(\bar{X}^{k}(t)))\mathrm{d}B^{0}_t
\end{aligned}
\end{align}
with the initial data $\bar{X}^{k}_0=\varsigma^k$. Theorem 1.3 ensures the existence and uniqueness of the solutions of these equations.
By the similar techniques as \cite[Vol-II, Proposition 2.11]{CD2018}, one obtains that for any $k\geq 1$,
 \begin{align}\label{eq4.2*}
\mathbb{P}^0(\mathcal{L}^1(\bar{X}^{k}(t))=\mathcal{L}^1(\bar{X}^{1}(t)), \forall t\geq 0)=1.
\end{align}
We now give the corresponding functional interacting particle systems
\begin{align}\label{eq4.1}
\mathrm{d}X^{k,N}(t)&=f(X^{k,N}_t,\mu_t^N)\mathrm{d}t+g(X^{k,N}_t,\mu_t^N)\mathrm{d}B^k_t
\nonumber\\&\quad
+g^0(X^{k,N}_t,\mu_t^N)\mathrm{d}B^{0}_t,~k=1,\cdots,N,
\end{align}
with the initial data $ {X}^{k, N}_0=\varsigma^k$,
where $\mu_t^N=\frac{1}{N}\sum_{k=1}^N\delta_{X^{k,N}(t)}$ is the empirical measure of the particles.
Under Assumption 1, by the same way as in \cite[Proposition 3.1]{WXZ2023}, we can  prove that the equation \eqref{eq4.1} has a unique strong solution $\mathbf{X}^N(t)=(X^{1,N}(t),\cdots,X^{N,N}(t))$ on $[-\tau,+\infty)$.
    Now, we give the result on conditional propagation of chaos.
\begin{thm}\label{T2}
 Let {\rm Assumption \ref{A1}} and {\rm \ref{A2}} hold and  $p>2$, $\varsigma\in \mathbb{L}_p(C) $. Then there exists a constant $K$  depending on $T,d,{p,q}~(2\leq q<p)$ such that
\begin{align}
\label{eq4.4}&\max_{1\leq k\leq N}\mathbb{E}\Big[\sup_{t\in[0,T]}|X^{k,N}(t)-\bar{X}^k(t)|^{q}\Big]\leq K\varepsilon_{N},
\\&
\label{eq4.5}\sup_{t\in[0,T]}\mathbb{E}\Big[\mathbb{W}_{q}^{q}\big(\mu_t^N,\mathcal{L}^1(\bar{X}^{k}(t))\big)\Big]\leq K\varepsilon_{N},
\end{align}
where $\bar{X}^k(t),~X^{k,N}(t)$ are solutions to equations \eqref{eq4.2} and \eqref{eq4.1}, 
respectively, and $\varepsilon_{N}$ is defined by \eqref{L**}.
\end{thm}
{\bf Proof.}
Fix $2\leq q<p$. Let $Z^{k,N}(t)=X^{k,N}(t)-\bar{X}^k(t)$ and $Z^{k,N}_t=X^{k,N}_t-\bar{X}^k_t$. By It\^o formula and Assumption \ref{A1},
we derive from equations \eqref{eq4.2} and \eqref{eq4.1} that for any $t\in[0,T]$,
\begin{align*}
&\mathbb{E}\Big[\sup_{s\in[0,t]}|Z^{k,N}(s)|^{q}\Big]
\\&\leq \mathbb{E}\Big[\int_0^t\!\frac{2L_1 q\!+\!L_2 q(q\!-\!1)}{2}|Z^{k,N}(s)|^{q-2}(\|Z^{k,N}_s\|^2_C
+\mathbb{W}_2^2(\mu^{N}_s,\mathcal{L}^1(\bar{X}^k(s))))\mathrm{d}s
\\&+ \mathbb{E}\Big[\sup_{0\leq s\leq t}\int_0^s\!\!\!q|Z^{k,N}(r)|^{q-2}\langle Z^{k,N}(r),g(X^{k,N}_r,\mu^{N}_r)-g(\bar{X}^k_r,\mathcal{L}^1(\bar{X}^{k}(r)))\rangle\mathrm{d}B_r\Big]
\\&+\mathbb{E}\Big[\sup_{0\leq s\leq t}\int_0^s\!\!\!q|Z^{k,N}(r)|^{q-2}
\langle Z^{k,N}(r),g^0(X^{k,N}_r,\mu^{N}_r)
-\!g^0(\bar{X}^k_r,\mathcal{L}^1(\bar{X}^{k}(r)))\rangle\mathrm{d}B^0_r\Big].
\end{align*}
By virtue of the Burkholder-Davis-Gundy inequality, one arrives at
\begin{align*}
&\mathbb{E}\Big[\sup_{s\in[0,t]}|Z^{k,N}(s)|^{q}\Big]\\
&
\leq \mathbb{E}\int_0^t\frac{2L_1 q+L_2 q(q-1)}{2}|Z^{k,N}(s)|^{q-2}(\|Z^{k,N}_s\|^2_C
+\mathbb{W}_2^2(\mu^{N}_s,\mathcal{L}^1(\bar{X}^k(s))))\mathrm{d}s
\\
&+4\sqrt{2}q\mathbb{E}\Big(\int_0^t|Z^{k,N}(s)|^{2q-2}|g(X^{k,N}_s,\mu^{k,N}_s)-g(\bar{X}^k_s,\mathcal{L}^1(\bar{X}^{k}(s)))|^2\mathrm{d}s\Big)^{\frac{1}{2}}
\\
&+4\sqrt{2}q\mathbb{E}\Big(\int_0^t|Z^{k,N}(s)|^{2q-2}|g^0(X^{k,N}_s,\mu^{k,N}_s)\!-\!g^0(\bar{X}^k_s,\mathcal{L}^1(\bar{X}^{k}(s)))|^2\mathrm{d}s\Big)^{\frac{1}{2}}\\
&
\leq \mathbb{E}\int_0^t\Big(\frac{2L_1 q +L_2 q(q-1)}{2}\Big)|Z^{k,N}(s)|^{q-2}
(\|Z^{k,N}_s\|^2_C+\mathbb{W}_2^2(\mu^{N}_s,\mathcal{L}^1(\bar{X}^k(s))))\mathrm{d}s
\\
&+\!8\sqrt{2L_2}q\mathbb{E}\Big(\int_0^t|Z^{k,N}(s)|^{2q-2}(\|Z^{k,N}_s\|_C^2+\mathbb{W}_2^2(\mu^{N}_s,\mathcal{L}^1(\bar{X}^k(s))))\mathrm{d}s\Big)^{\frac{1}{2}}
.
\end{align*}
Using the Young inequality, Assumption \ref{A1} and the discrete H\"older inequality yields
\begin{align*}
&\mathbb{E}\Big[\sup_{s\in[0,t]}|Z^{k,N}(s)|^{q}\Big]
\\&\leq \mathbb{E}\int_0^t\frac{2L_1 q+L_2 q(q-1)}{2}|Z^{k,N}(s)|^{q-2}(\|Z^{k,N}_s\|^2_C+\mathbb{W}_2^2(\mu^{N}_s,\mathcal{L}^1(\bar{X}^k(s))))\mathrm{d}s
\\
&+8\sqrt{2L_2}q\mathbb{E}\Big(\sup_{s\in[0,t]}|Z^{k,N}(s)|^{q}\int_0^t\!\!|Z^{k,N}(s)|^{q-2}
(\|Z^{k,N}_s\|_C^2+\mathbb{W}_2^2(\mu^{N}_s,\mathcal{L}^1(\bar{X}^k(s))))\mathrm{d}s\Big)^{\frac{1}{2}}
\end{align*}
\begin{align*}
&
\leq\frac{1}{2}\mathbb{E}\Big[\sup_{s\in[0,t]}|Z^{k,N}(s)|^{q}\Big]+\mathbb{E}\int_0^t\Big(\frac{2L_1 q+L_2 q(257q-1)}{2}\Big)|Z^{k,N}(s)|^{q-2}
\\
&\quad \times(\|Z^{k,N}_s\|^2_C+\mathbb{W}_2^2(\mu^{N}_s,\mathcal{L}^1(\bar{X}^k(s))))\mathrm{d}s
\\
&\leq\frac{1}{2}\mathbb{E}\Big[\sup_{s\in[0,t]}|Z^{k,N}(s)|^{q}\Big]+\big(2L_1+L_2(257 q-1)\big)\big(q-1\big)\mathbb{E}\int_0^t
\|Z^{k,N}_s\|_C^{q}\mathrm{d}s
\\&\quad+ \big(2L_1+L_2(257 q-1)\big)\mathbb{E}\int_0^t\mathbb{W}_2^q(\mu^{N}_s,\mathcal{L}^1(\bar{X}^k(s)))\mathrm{d}s.
\end{align*}
Therefore,
\begin{align*}
\begin{aligned}
\mathbb{E}\Big[\sup_{s\in[0,t]}|Z^{k,N}(s)|^{q}\Big]
\leq K_q\mathbb{E}\int_0^t\sup_{r\in[0,s]}
\|Z^{k,N}_r\|^q_C\mathrm{d}s+K_q\mathbb{E}\int_0^t\mathbb{W}_2^q(\mu^{N}_s,\mathcal{L}^1(\bar{X}^k(s)))\mathrm{d}s,
\end{aligned}
\end{align*}
where $K_q=2\big(2L_1+L_2(257 q-1)\big)\big(q-1\big)$.
As the initial values of $X^{k,N}_s$ and $\bar{X}^k_s$ are the same, we derive from the Fubini theorem that
\begin{align}\label{eq8}
\begin{aligned}
&\mathbb{E}\Big[\sup_{s\in[0,t]}|Z^{k,N}(s)|^{q}\Big]
\\&\leq K_q\int_0^t\mathbb{E}\Big[\sup_{r\in[0,s]}
|Z^{k,N}(r)|^q\Big]\mathrm{d}s+{K}_q\int_0^t\mathbb{E}\mathbb{W}_2^q(\mu^{N}_s,\mathcal{L}^1(\bar{X}^k(s)))\mathrm{d}s.
\end{aligned}
\end{align}
Using the H\"older inequality, the triangle inequality and the discrete H\"older inequality derives that
\begin{align}\label{eq10}
\begin{aligned}
\mathbb{E}\mathbb{W}_{2}^{q}(\mu_t^N,\mathcal{L}^1(\bar{X}^{k}(t)))
&\leq\mathbb{E}\mathbb{W}_{q}^{q}(\mu_t^N,\mathcal{L}^1(\bar{X}^{k}(t)))
\\&\leq \mathbb{E}\big(\mathbb{W}_{q}(\mu_t^N,\nu_t^N)+\mathbb{W}_{q}(\nu_t^N,\mathcal{L}^1(\bar{X}^{k}(t)))\big)^q
\\&\leq 2^{q-1}\mathbb{E}\mathbb{W}_{q}^{q}(\mu_t^N,\nu_t^N)+2^{q-1}\mathbb{E}\mathbb{W}_{q}^{q}(\nu_t^N,\mathcal{L}^1(\bar{X}^{k}(t))),
\end{aligned}
\end{align}
 where $\nu_t^N\!=\!\frac{1}{N}\sum_{k=1}^N\delta_{\bar{X}^{k}(t)}$.
Owing to the symmetrical structure of \eqref{eq4.1}, and noticing that the distributions of  $\bar{X}^k(t),~k=1,\cdots, N$ are identical, we compute
\begin{align}\label{eq11}
\begin{aligned}
\mathbb{E}\mathbb{W}_{q}^{q}(\mu_t^N,\nu_t^N)&=\mathbb{E}\mathbb{W}_{q}^{q}\Big(\frac{1}{N}\sum_{k=1}^N\delta_{{X}^{k,N}(t)},\frac{1}{N}\sum_{k=1}^N\delta_{\bar{X}^{k}(t)}\Big)
\\&\leq\mathbb{E}\Big[\frac{1}{N}\sum_{k=1}^N|Z^{k,N}(t)|^{q}\Big]
\\&\leq \frac{1}{N}\sum_{k=1}^N\mathbb{E}|Z^{k,N}(t)|^{q}
\\&\leq\mathbb{E}|Z^{k,N}(t)|^{q}.
\end{aligned}
\end{align}
It follows from \eqref{eq4.2*} that there exists a set $\tilde{\Omega}^0$ satisfying $\mathbb{P}^0(\tilde{\Omega}^0)=1$ such that for any $\omega^0\in\tilde{\Omega}^0$,
$$\mathcal{L}^1(\bar{X}^{k}(t))(\omega^0)=\mathcal{L}^1(\bar{X}^{1}(t))(\omega^0),\quad\forall t\in[0,T],$$
which means that for any $\omega^0\in\tilde{\Omega}^0$ and $t\in[0,T]$, $ \bar{X}^{k}(t,\omega^0,\cdot),   k=1,\cdots,   N, $ are independent and identically distributed random variables with distribution $\mathcal{L}^1(\bar{X}^{1}(t))(\omega^0)$. Then one deduces from Lemma \ref{L*} that for any $\omega^0\in\tilde{\Omega}^0$, for any $
0\leq t\leq T$,
\begin{align*}
\begin{aligned}
&\mathbb{E}^1 \mathbb{W}_{q}^{q}\big(\nu_t^N(\omega^0),\mathcal{L}^1(\bar{X}^{k}(t))(\omega^0)\!\big)
\\&=\mathbb{E}^1 \mathbb{W}_{q}^{q}\Big(\frac{1}{N}\sum_{k=1}^N\delta_{\bar{X}^{k}(t,\omega^0,\cdot)},\mathcal{L}^1(\bar{X}^{k}(t))(\omega^0)\!\Big)
\\
&\leq c\varepsilon_{N}\mathbb{W}_{p}^{q}(\mathcal{L}^1(\bar{X}^{k}(t))(\omega^0),\delta_0)
\leq c\varepsilon_{N}(\mathbb{E}^1|\bar{X}^1(t)|^p)^{\frac{q}{p}},
\end{aligned}
\end{align*}
where $c$ is a positive constant depending on $d,p,q$, and $\varepsilon_{N}$ is defined by \eqref{L**}. Then using the property of the expectation and the H\"older inequality, one arrives at
\begin{align}\label{eq12*}
\mathbb{E}\mathbb{W}_{q}^{q}\big(\nu_t^N,\mathcal{L}^1(\bar{X}^{k}(t))\big)
&=\mathbb{E}^0[\mathbb{E}^1 \mathbb{W}_{q}^{q}\big(\nu_t^N,\mathcal{L}^1(\bar{X}^{k}(t))\big)]
\nonumber\\&\leq c\varepsilon_{N}\mathbb{E}^0[(\mathbb{E}^1|\bar{X}^1(t)|^p)^{\frac{q}{p}}]
\nonumber\\&\leq c\varepsilon_{N}(\mathbb{E}^0[\mathbb{E}^1|\bar{X}^1(t)|^p])^{\frac{q}{p}}
\nonumber\\&\leq c\varepsilon_{N}(\mathbb{E}|\bar{X}^1(t)|^p)^{\frac{q}{p}}.
\end{align}
Inserting \eqref{eq11} and \eqref{eq12*} into \eqref{eq10} shows that
\begin{align}\label{eq13}
\mathbb{E}\mathbb{W}_{2}^{q}(\mu_t^N,\mathcal{L}^1(\bar{X}^{k}(t)))
\leq 2^{q-1}\mathbb{E}|Z^{k,N}(t)|^{q}
+ 2^{q-1}c\varepsilon_{N}(\mathbb{E}|\bar{X}^1(t)|^p)^{\frac{q}{p}}.
\end{align}
Combining \eqref{eq8} and \eqref{eq13}, we obtain that
\begin{align*}
\mathbb{E}\Big[\sup_{s\in[0,t]}|Z^{k,N}(s)|^{q}\Big]
&\leq K_1\int_0^t\mathbb{E}\Big[\sup_{r\in[0,s]}
|Z^{k,N}(r)|^q\Big]\mathrm{d}s+K_1\varepsilon_{N}\int_0^t\mathbb{E}|\bar{X}^1(s)|^p\mathrm{d}s
\\&
\leq K_1
\int_0^t\mathbb{E}\Big[\sup_{r\in[0,s]}|Z^{k,N}(r)|^{q}\Big]\mathrm{d} s
+K_1\varepsilon_{N}T\Big(\mathbb{E}\Big[\sup_{t\in[0,T]}|\bar{X}^1(t)|^p\Big]\Big)^{\frac{q}{p}},
\end{align*}
where $K_1=\max\{(1+2^{q-1})K_q,c2^{q-1}K_q\}$.
Using the Gronwall inequality and Theorem \ref{T1} gives that
\begin{align}\label{eq4.6*}
\begin{aligned}
&\mathbb{E}\Big[\sup_{t\in[0,T]}|Z^{k,N}(t)|^{q}\Big]
\leq K_1Te^{K_1T}\varepsilon_{N}\Big(\mathbb{E}\Big[\sup_{t\in[0,T]}|\bar{X}^1(t)|^p\Big]\Big)^{\frac{q}{p}}
\leq K \varepsilon_{N}.
\end{aligned}
\end{align}
where $K=K_1Te^{K_1T}\Big(\mathbb{E}\Big[\sup_{s\in[0,T]}|\bar{X}^1(t)|^p\Big]\Big)^{{q}/{p}}$.
Taking the maximum value w.r.t. $k$ on both sides of \eqref{eq4.6*} gives the assertion \eqref{eq4.4}. Now we turn to prove the assertion \eqref{eq4.5}.
 Ultilizing \eqref{eq13}, \eqref{eq4.6*} and Theorem \ref{T1}, one obtains that
\begin{align*}
&\sup_{t\in[0,T]}\mathbb{E}\mathbb{W}_{q}^{q}(\mu_t^N,\mathcal{L}^1(\bar{X}^{k}(t)))
\\&\leq 2^{q-1}\sup_{t\in[0,T]}\mathbb{E}|Z^{k,N}(t)|^{q}+2^{q-1}\varepsilon_{N}\sup_{t\in[0,T]}(\mathbb{E}|\bar{X}^1(t)|^p)^{\frac{q}{p}}
\\&\leq 2^{q-1}\mathbb{E}\Big[\sup_{t\in[0,T]}|Z^{k,N}(t)|^{q}\Big]+2^{q-1}\varepsilon_{N}\Big(\mathbb{E}\Big[\sup_{t\in[0,T]}|\bar{X}^1(t)|^p\Big]\Big)^{\frac{q}{p}}
\\&\leq K\varepsilon_{N},
\end{align*}
where $K=2^{q-1}(K_1Te^{K_1T}+1)\Big(\mathbb{E}\Big[\sup_{t\in[0,T]}|\bar{X}^1(t)|^p\Big]\Big)^{{q}/{p}}$.
Then the desired assertions hold.
\qed

Utilizing the conditional propagation of chaos, we now investigate the stability equivalence between the equation  \eqref{eq4.2} and the corresponding functional particle systems with common noise \eqref{eq4.1}.
\begin{thm}\label{T3}
 Let {\rm Assumption \ref{A1}} and {\rm \ref{A2}} hold and  $p>2$, $\varsigma\in \mathbb{L}_p(C) $. Then the solution of \eqref{eq4.2} is exponentially stable in the $q$th moment ($2\leq q< p$), i.e. there exists a positive constant $\kappa_1$ such that for any $k>1$,
 \begin{align}\label{T3*}
 \limsup _{t \rightarrow \infty} \frac{1}{t} \log \left(\mathbb{E}|\bar{X}^k(t)|^q\right) \leq-\kappa_1
 \end{align}
if and only if the solution of \eqref{eq4.1} is exponentially stable in the $q$th  moment ($2\leq q< p$), i.e. there exists a positive constant $\kappa_2$ such that for any $k>1$,
 \begin{align}\label{T3**}
 \limsup _{t \rightarrow \infty}\lim_{N\rightarrow\infty} \frac{1}{t} \log \left(\mathbb{E}|{X}^{k,N}(t)|^q\right) \leq-\kappa_2.
 \end{align}
\end{thm}
{\bf Proof.}
We first prove that \eqref{T3*} implies \eqref{T3**}. Making use of the discrete H\"older inequality gives for any $k\geq1$
\begin{align*}
\mathbb{E}|{X}^{k,N}(t)|^q\leq 2^{q-1}\mathbb{E}|{X}^{k,N}(t)-\bar{X}^{k}(t)|^q+2^{q-1}\mathbb{E}|\bar{X}^{k}(t)|^q, \quad t\geq 0.
\end{align*}
Taking limitation with respect to $N$ on both sides and using Theorem \ref{T2}, one arrives at that for any $t\geq0$,
\begin{align*}
\lim_{N\rightarrow\infty}\mathbb{E}|{X}^{k,N}(t)|^q&\leq2^{q-1}\lim_{N\rightarrow\infty}\mathbb{E}|{X}^{k,N}(t)-\bar{X}^{k}(t)|^q+2^{q-1}\mathbb{E}|\bar{X}^{k}(t)|^q
\\&\leq 2^{q-1}\mathbb{E}|\bar{X}^{k}(t)|^q.
\end{align*}
Then it follows from taking the logarithm, dividing $t$ and taking limitation with respect to $t$ on both sides that
\begin{align*}
\lim_{t\rightarrow\infty}\frac{1}{t}\log\big(\lim_{N\rightarrow\infty}\mathbb{E}|{X}^{k,N}(t)|^q\big)\leq \lim_{t\rightarrow\infty}\frac{1}{t}(q-1)\log(2)+\lim_{t\rightarrow\infty}\frac{1}{t}\log(\mathbb{E}|\bar{X}^{k}(t)|^q)
\leq \kappa_1.
\end{align*}
The continuity of the logarithm function gives \eqref{T3**}.
By the similar procedure, we obtain that \eqref{T3**} implies \eqref{T3*}. The proof is therefore complete.
\qed
\section{The stability}\label{Sect.4}
This section focuses on the stability of the MV-SFDEs with common noise \eqref{1.1}. We first give the It\^o formula involved with state and measure for the equation \eqref{1.1}. Then we establish the Razumikhin theorem for the exponential stability of MV-SFDEs with common noise \eqref{1.1} using the It\^o formula involved with state and measure. We now impose another Assumption. 
\begin{assp}
 \label{A3}
 \rm Assume that $f({0},\delta_{0})=g(0,\delta_{0})=g^0({0},\delta_{0})=0$.
  There exist positive constants ${L}_5$ and $l\geq2$, such that for any $\phi\in C$ and any $\mu\in\mathcal{P}_{2}(\mathbb{R}^d)$,
\begin{align*}
&|{f}(\phi, \mu)|^2
 \leq {L}_5(1+\|\phi\|_C^l+\mathbb{W}^2_2(\mu,\delta_0)).
\end{align*}
\end{assp}
For convenience, we give the following notations. For any $p\geq2,t\geq0$, denote by $\mathbb{L}_{t,p}(C)$ the family of all $\mathcal{F}_t$-measurable $C$-valued random variable $\phi$ satisfying $\mathbb{E}\|\phi\|_C^p<\infty$. Denote by $\mathbb{L}_{t}^0( \mathcal{P}_p )$ the family of all $\mathcal{F}^0_t$-measurable $\mathcal{P}_p (\mathbb{R}^d)$-valued random variable $\mu$ on $\Omega^0$ such that $\mathbb{E}^0\mathbb{W}_p^p(\mu,\delta_0)<\infty$. For each operator $V\in\tilde{C}^{2,2,1}$, define a  random operator $ LV$
by
\begin{align} \label{lv}
LV(\phi,\mu,t)
& =\partial_t V(\phi(0),\mu,t)
+\partial_x V(\phi(0),\mu,t) {f}(\phi,\mu)
\nonumber\\&\quad+\frac{1}{2} \operatorname{tr}\big( \partial_{x x}V(\phi(0),\mu,t){g}(\phi,\mu)
({g}(\phi,\mu))^{\prime} \big)
\nonumber\\
&\quad+\frac{1}{2} \operatorname{tr}\!\big( \partial_{x x}V(\phi(0),\mu,t){g}^0(\phi,\mu)
({g}^0(\phi,\mu))^{\prime} \big)
\nonumber\\
&\quad+\mathbb{E}_1^1\Big[\big(\partial_\mu V(\phi(0),\mu,t)(\phi_1(0))\big)^{\prime} {{{f}}}(\phi_1,\mu_1)
\nonumber\\
&
\quad+  \operatorname{tr}\big(\partial_x \partial_\mu V(\phi(0),\mu,t)(\phi_1(0)) {g}^0(\phi,\mu)
({g}^0(\phi_1,\mu_1))^{\prime}\big)
\nonumber\\
&\quad+\frac{1}{2} \operatorname{tr}\big(\partial_y \partial_\mu V(\phi(0),\mu,t)(\phi_1(0)) {{g}(\phi_1,\mu_1)}({{g}(\phi_1,\mu_1)})^{\prime}\big)
\nonumber\\
&
\quad+\frac{1}{2} \operatorname{tr}\big(\partial_y \partial_\mu V(\phi(0),\mu,t)(\phi_1(0)) {{g}}^0(\phi_1,\mu_1)
({{g}}^0(\phi_1,\mu_1))^{\prime}\big)
\nonumber\\
&
\quad+\frac{1}{2} \mathbb{E}_2^1 \Big[\operatorname{tr}\big(\partial_\mu^2 V(\phi(0),\mu,t)(\phi_1(0), {\phi_2(0)}){{g}}^0(\phi_1,\mu_1)({{g}}^0(\phi_2,\mu_2))^{\prime}\big)
\Big]\Big],
\end{align}
 where $\phi_i$, $\mu_i$ are copies of $\phi$, $\mu$ on $\Omega_i~(i=1,2)$, for   $(\phi,\mu,t)\in \mathbb{L}_2(C)\times \mathcal{P}_2(\mathbb{R}^d)\times [0,\infty)$ such that the right conditional expectations  exist.

 \begin{lem}\label{le5.1}
 Let {\rm Assumption \ref{A1}, \ref{A2} and \ref{A3}} hold and $V\in\tilde{C}^{2,2,1}$. Then the solution $X(t)$ to the MV-SFDE with common noise \eqref{1.1} with initial data $\varsigma\in \mathbb{L}_p(C)$ ($p\geq4 \vee l$) satisfies for any $t\geq0$
\begin{align*}
\mathbb{E} V(X(t),\mathcal{L}^1(X(t)),t)
&=\mathbb{E} V(\varsigma(0),\mathcal{L}^{1  }(\varsigma(0)),0)
+\mathbb{E}\int_0^t{L}V(X_s,\mathcal{L}^1(X(s)),s)\mathrm{d}s.
\end{align*}
\end{lem}
{\bf Proof.}
Let $X(t)$ be the solution to the equation \eqref{1.1}. Due to \cite[Vol-II, Theorem 4.17]{CD2018}, it is sufficient to prove that for any $T\geq 0$,
\begin{align}\label{eq5.3*}
&\mathbb{E}\int_0^T\Big[|{f}(X_t,\mathcal{L}^1(X(t)))|^2
+|{g}(X_t,\mathcal{L}^1(X(t)))|^4
+|{g}^0(X_t,\mathcal{L}^1(X(t)))|^4\Big]\mathrm{d}t<\infty.
\end{align}
We notice 
from Assumption \ref{A3} and Theorem \ref{T1} that for any $T\geq 0$,
 \begin{align*}
\mathbb{E}\int_0^{T}|{f}(X_t,\mathcal{L}^1(X(t)))|^2\mathrm{d}t
&\leq {L}_5\mathbb{E}\int_0^{T}(1+\|X_t\|_{C}^l+\mathbb{W}^2_2(\mathcal{L}^1(X(t)),\delta_0))\mathrm{d}t
\nonumber\\&\leq {L}_5\mathbb{E}\int_0^{T}(1+\|X_t\|_{C}^l+\mathbb{W}^2_2(\mathcal{L}^1(X(t)),\delta_0))\mathrm{d}t
\nonumber\\
&\leq{L}_5+2{L}_5T\mathbb{E}\Big[\sup_{-\tau\leq t\leq T}|X(t)|^l\Big]<\infty.
\end{align*}
By \eqref{A11}, we have
\begin{align*}
&\mathbb{E}\int_0^T(|{g}(X_t, \mathcal{L}^1(X(t)))|^4
+|\bar{g}^0(X_t,\mathcal{L}^1(X(t)))|^4) \mathrm{d} t
\nonumber\\&
\leq 3{L}_4^2\mathbb{E}\int_0^T(1+\|X_t\|_C^4+\mathbb{E}^1|X(t)|^4)\mathrm{d}t
\nonumber\\&
\leq 3{L}_4^2T+6{L}_4^2T\mathbb{E}\Big[\sup_{-\tau\leq t\leq T}|X(t)|^4\Big]
<\infty.
\end{align*}
 This implies that \eqref{eq5.3*} holds. Then for $V\in\tilde{C}^{2,2,1}$, using \cite[Vol-II, Theorem 4.17]{CD2018} and taking expectation on both sides, one has for any $t\geq0$, $\mathbb{P}$-a.s.
\begin{align*}
&\E V(X(t),\mathcal{L}^1(X(t))),t)
\nonumber\\&=\E V(\varsigma(0),\mathcal{L}^1(\varsigma(0)),0)+\E\int_0^t\Big( \partial_t V(X(s),\mathcal{L}^1(X(s)),s)
\nonumber\\
&
+\partial_x V(X(s),\mathcal{L}^1(X(s)),s){f}(X_s,\mathcal{L}^1(X(s)))
\nonumber\\&+\frac{1}{2} \operatorname{tr}\big(\partial_{x x}V(X(s),\mathcal{L}^1(X(s)),s){g}(X_s,\mathcal{L}^1(X(s)))
({g}(X_s,\mathcal{L}^1(X(s))))^{\prime}\big)
\nonumber\\&+\frac{1}{2} \operatorname{tr}\big(\partial_{x x}V(X(s),\mathcal{L}^1(X(s)),s){g}^0(X_s,\mathcal{L}^1(X(s)))
({g}^0(X_s,\mathcal{L}^1(X(s))))^{\prime}\big)
\nonumber\\
&+\mathbb{E}_1^1\Big[\big(\partial_\mu V(X(s),\mathcal{L}^1(X(s)),s)(X_1(s))\big)^{\prime} {f}(X_{s,1},\mathcal{L}^1(X_1(s)))
\nonumber\\&+ \operatorname{tr}\big(\partial_x \partial_\mu V(X(s),\mathcal{L}^1(X(s)),s)(X_1(s)) {g}^0(X_{s},\mathcal{L}^1(X(s)))
({g}^0(X_{s,1},\mathcal{L}^1(X_1(s))))^{\prime}\big)
\nonumber\\
&+\frac{1}{2} \operatorname{tr}\big(\partial_y \partial_\mu V(X(s),\mathcal{L}^1(X(s)),s)(X_1(s)) {{g}(X_{s,1},\mathcal{L}^1(X_1(s)))}({g}(X_{s,1},\mathcal{L}^1(X_1(s))))^{\prime}\big)
\nonumber\\
&+ \frac{1}{2} \operatorname{tr}\big(\partial_y \partial_\mu V(X(s),\mathcal{L}^1(X(s)),s)(X_1(s)){g}^0(X_{s,1},\mathcal{L}^1(X_1(s)))
({g}^0(X_{s,1},\mathcal{L}^1(X_1(s))))^{\prime}\big)
\nonumber\\
&+\frac{1}{2} \mathbb{E}_2^1\big[\operatorname{tr}\big(\partial_\mu^2 V(X(s),\mathcal{L}^1(X(s)),s)(X_1(s), {X_2}(s)){g}^0(X_{s,1},\mathcal{L}^1(X_1(s)))
\nonumber\\
&\quad\quad\quad\quad\quad\quad\times({g}^0(X_{s,2},\mathcal{L}^1(X_2(s))))^{\prime}\big)\big]\Big]\Big)\mathrm{d} s,
\end{align*}
where $X_1(t)$, $X_2(t)$ are copies of $X(t)$ on probability spaces $\Omega_1$ and $\Omega_2$, respectively, and $\mathcal{L}^1(X_1(t))$, $\mathcal{L}^1(X_2(t))$ are copies of $\mathcal{L}^1(X(t))$ on probability spaces $\Omega_1$ and $\Omega_2$, respectively.
Then let $X_{t,1}$, $X_{t,2}$ be the segment of $X_1(t)$ and $X_2(t)$, respectively. It is obvious that $X_{t,1}$, $X_{t,2}$ are the copies of the segment $X_t$ of $X(t)$ on probability spaces $\Omega_1$ and $\Omega_2$. Then the desired assertion holds.\qed

We proceed to establish the Razumikhin theorem for the exponential stability of the MV-SFDEs with common noise \eqref{1.1}.

\begin{thm}\label{th3.2}
Let {\rm Assumption \ref{A1}, \ref{A2}} and {\rm \ref{A3}} hold. Assume that there exists a function $V\in \tilde{C}^{2,2,1}$ and positive constants $\lambda$, $q$, $c_1,$ $c_2$, $c_3$, $c_4$, and $\alpha>1$, such that
\begin{align}\label{eq5.5}
\begin{aligned}
c_1|x|^q+c_2\mathbb{W}^q_q(\mu,\delta_0)&\leq V(x,\mu,t)
\leq c_3 |x|^q+c_4\mathbb{W}^q_q(\mu,\delta_0),
\end{aligned}
\end{align}
for any $x\in\mathbb{R}^d$, $\mu\in \mathcal{P}_q(\mathbb{R}^d)$, and $t\in[-\tau,\infty)$.
Moreover,
\begin{align}\label{eq5.6}
\mathbb{E}LV(\phi,\mathcal{L}^1(\phi(0)),t)&\leq -\lambda \mathbb{E} V(\phi(0),\mathcal{L}^1(\phi(0)),t)
,
\end{align}
for any $t\geq0$ and those $\phi\in \mathbb{L}_{t,q}(C)$ satisfying
\begin{align*}
\mathbb{E}V(\phi(\theta),\mathcal{L}^1(\phi(\theta)),t+\theta)<\alpha\mathbb{E}V(\phi(0),\mathcal{L}^1(\phi(0)),t), \quad \mathrm{on} ~-\tau\leq \theta\leq0.
\end{align*}
Then the solution $X(t)$ to \eqref{1.1} with the initial value $\varsigma\in\mathbb{L}_{p}(C)$($p\geq4 \vee l \vee q$) satisfies
\begin{align}\label{eq5.8}
\mathbb{E}|X(t)|^q\leq \frac{c_3+c_4}{c_1+c_2}\mathbb{E}\|\varsigma\|^qe^{-\kappa t},\quad \mathrm{on}~t\geq0,
\end{align}
where $\kappa=\min\{\lambda,\log(\alpha)/\tau\}$.
\end{thm}
{\bf Proof.}
Let $X(t)$ be the solution to the MV-SFDE with common noise \eqref{1.1} with initial value $\varsigma\in\mathbb{L}_{p}(C)$.  For any $\varepsilon\in(0,\kappa)$, let $\bar{\kappa}=\kappa-\varepsilon$. Define
\begin{align}\label{eq5.9}
I(t)=\max_{\theta\in[-\tau,0]}\big[e^{\bar{\kappa}(t+\theta)}\mathbb{E}V(X(t+\theta),\mathcal{L}^1(X(t+\theta)),t+\theta)\big], ~t\geq 0.
\end{align}
Thanks to Theorem \ref{T1}, we have for any $T\geq0$, $\mathbb{E}\Big[\sup_{t\in[0,T]}|X(t)|^p\Big]<\infty$. Note that $X(t)$ is continuous with respect to $t$. Then it follows from Lemma \ref{l2} that $\mathcal{L}^1(X(t))$ is continuous with respect to $t$. The continuity with respect to $t$ of $X(t)$, $\mathcal{L}^1(X(t))$ and $V(x,\mu,t)$ gives that $\mathbb{E}V(X(t),\mathcal{L}^1(X(t)),t)$ is continuous with respect to $t$. Then we know that $I(t)$ is well-defined and continuous with respect to $t$.Define
\begin{align}\label{eq5.91}
D^+I(t):=\limsup_{h\rightarrow 0^+}\frac{I(t+h)-I(t)}{h},\quad \mathrm{for~ any}~t\geq0.
\end{align}
If \begin{align}\label{eq5.912}
D^+I(t)\leq 0, \quad \mathrm{for~ any}~t\geq0,
\end{align}
we obtain that
\begin{align}\label{eq5.92}
I(t)\leq I(0),
\end{align} for any $t\geq0$. Using \eqref{eq5.92}, the definition of $I(t)$ in \eqref{eq5.9} and the inequality \eqref{eq5.5} shows that
\begin{align}\label{eq5.10}
\begin{aligned}
&\max_{\theta\in[-\tau,0]}\big[e^{\bar{\kappa}(t+\theta)}\mathbb{E}V(X(t+\theta),\mathcal{L}^1(X(t+\theta)),t+\theta)\big]
\\&\leq \max_{\theta\in[-\tau,0]}\mathbb{E}V(X(\theta),\mathcal{L}^1(X(\theta)),\theta)
\\&\leq (c_3+c_4)\max_{\theta\in[-\tau,0]}\mathbb{E}|\varsigma(\theta)|^q
\\&\leq (c_3+c_4)\mathbb{E}\|\varsigma\|_C^q.
\end{aligned}
\end{align}
Using \eqref{eq5.5} and \eqref{eq5.10}, one obtains that
\begin{align}\label{eq5.11}
\begin{aligned}
(c_1+c_2)e^{\bar{\kappa}t}\mathbb{E}|X(t)|^q
&\leq e^{\bar{\kappa}t}\mathbb{E}V(X(t),\mathcal{L}^1(X(t)),t)
\\&\leq\max_{\theta\in[-\tau,0]}\big[e^{\bar{\kappa}(t+\theta)}\mathbb{E}V(X(t+\theta),\mathcal{L}^1(X(t+\theta)),t+\theta)\big]
\\&\leq (c_3+c_4)\mathbb{E}\|\varsigma\|_C^q.
\end{aligned}
\end{align}
Dividing $(c_1+c_2)e^{\bar{\kappa}t}$ on both side of \eqref{eq5.11} yields that
\begin{align*}
\mathbb{E}|X(t)|^q
\leq \frac{c_3+c_4}{c_1+c_2}\mathbb{E}\|\varsigma\|_C^q e^{-\bar{\kappa}t}.
\end{align*}
Owing to the arbitrariness of  $\varepsilon$,
  we obtain the desired assertion \eqref{eq5.8}. It is sufficient to prove \eqref{eq5.912} holds. 
 For any $t\geq 0$ (fixed at this time), define
\begin{align*} 
\tilde{\theta}=\max\{\theta\in[-\tau,0]:e^{\bar{\kappa}(t+\theta)}\mathbb{E}V(X(t+\theta),\mathcal{L}^1(X(t+\theta)),t+\theta)=I(t)\}.
\end{align*}
Then it is easy to see that $\tilde{\theta}$ is well-defined, $\tilde{\theta}\in[-\tau,0]$ and \begin{align*}
I(t)=e^{\bar{\kappa}(t+\tilde{\theta})}\mathbb{E}V(X(t+\tilde{\theta}),\mathcal{L}^1(X(t+\tilde{\theta})),t+\tilde{\theta}).
\end{align*}
We divide the proof into two cases.\\
\underline{Case 1} If $\tilde{\theta}<0$, then by the definition of $I(t)$, we have
\begin{align*}
e^{\bar{\kappa}(t+\theta)}\mathbb{E}V(X(t+\theta),\mathcal{L}^1(X(t+\theta)),t+\theta)< e^{\bar{\kappa}(t+\tilde{\theta})}\mathbb{E}V(X(t+\tilde{\theta}),\mathcal{L}^1(X(t+\tilde{\theta})),t+\tilde{\theta}),
\end{align*}
for any $\theta\in(\tilde{\theta},0]$. In particular, setting $\theta=0$, we have
\begin{align*}
e^{\bar{\kappa}t}\mathbb{E}V(X(t),\mathcal{L}^1(X(t)),t)< e^{\bar{\kappa}(t+\tilde{\theta})}\mathbb{E}V(X(t+\tilde{\theta}),\mathcal{L}^1(X(t+\tilde{\theta})),t+\tilde{\theta}),
\end{align*}
Since $e^{\bar{\kappa}t}\mathbb{E}V(X(t),\mathcal{L}^1(X(t)),t)$ is continuous with respect to $t$, then we can find a $h_0>0$ sufficient small such that for all $h\in[0,h_0]$ such that
\begin{align*}
e^{\bar{\kappa}(t+h)}\mathbb{E}V(X(t+h),\mathcal{L}^1(X(t+h)),t+h)\leq e^{\bar{\kappa}(t+\tilde{\theta})}\mathbb{E}V(X(t+\tilde{\theta}),\mathcal{L}^1(X(t+\tilde{\theta})),t+\tilde{\theta}),
\end{align*}
Thus one has
\begin{align*}
I(t+h)\leq I(t), ~h\in[0,h_0].
\end{align*}
Then this together with the definition of $D^+I(t)$ in \eqref{eq5.91} implies that $D^+I(t)\leq 0$. \\
\underline{Case 2} If $\tilde{\theta}=0$, then we arrive at that for any $\theta\in[-\tau,0]$
\begin{align}\label{eq5.12}
e^{\bar{\kappa}(t+\theta)}\mathbb{E}V(X(t+\theta),\mathcal{L}^1(X(t+\theta)),t+\theta)\leq e^{\bar{\kappa}t}\mathbb{E}V(X(t),\mathcal{L}^1(X(t)),t).
\end{align}
 Therefore, dividing $e^{\bar{\kappa}(t+\theta)}$ on both sides of \eqref{eq5.12} gives that for any $\theta\in[-\tau,0]$
\begin{align}\label{eq5.13}
\begin{aligned}
&\mathbb{E}V(X(t+\theta),\mathcal{L}^1(X(t+\theta)),t+\theta)\\&\leq e^{-\bar{\kappa}\theta}\mathbb{E}V(X(t),\mathcal{L}^1(X(t)),t)
\\&\leq e^{\bar{\kappa}\tau}\mathbb{E}V(X(t),\mathcal{L}^1(X(t)),t).
\end{aligned}\end{align}
 By virtue of \eqref{eq5.5}, one derives  that for any $\theta\in[-\tau,0]$
\begin{align}\label{eq5.14}
0\leq(c_1+c_2)\mathbb{E}|X(t+\theta)|^q\leq\mathbb{E}V(X(t+\theta),\mathcal{L}^1(X(t+\theta)),t+\theta ).
\end{align}
If $\mathbb{E}V(X(t),\mathcal{L}^1(X(t)),t)=0$, then by \eqref{eq5.13} and \eqref{eq5.14}, one has
\begin{align*}
\mathbb{E}|X(t+\theta)|^q&\leq \frac{1}{c_1+c_2}\mathbb{E}V(X(t+\theta),\mathcal{L}^1(X(t+\theta)),t+\theta )\\&\leq\frac{1}{c_1+c_2}e^{\bar{\kappa}\tau}\mathbb{E}V(X(t),\mathcal{L}^1(X(t)),t)=0,
\end{align*} for any $\theta\in[-\tau,0]$, which implies that $X(t+\theta)=0$ a.s. and $\mathcal{L}^1(X(t+\theta))=\delta_0$ a.s. for any $\theta\in[-\tau,0]$. Due to the uniqueness of the trivial solution $X(t)=0$, we derive that
$X(t+h)=0$ a.s. and $\mathcal{L}^1(X(t+h))=\delta_0$ a.s. for all $h>0$. Hence $D^+I(t)=0$. On the other hand, if $\mathbb{E}V(X(t),\mathcal{L}^1(X(t)),t)>0$, recalling the definition of $\kappa$, we know that $e^{\kappa \tau}\leq \alpha$, which together with \eqref{eq5.13} implies that
\begin{align*}
\mathbb{E}V(X(t+\theta),\mathcal{L}^1(X(t+\theta)),t+\theta)< \alpha\mathbb{E}V(X(t),\mathcal{L}^1(X(t)),t),
\end{align*}
for any $\theta\in[-\tau,0]$.  Then we derive from \eqref{eq5.6} that
\begin{align*}
\mathbb{E}{L}V(X_t,\mathcal{L}^1(X(t)),t)\leq -\lambda\mathbb{E}V(X(t),\mathcal{L}^1(X(t)),t).
\end{align*}
By the definition of $\kappa$, one has $\kappa\leq\lambda$.
Then the definition of $\bar{\kappa}$ shows that
\begin{align*}
\begin{aligned}
&\bar{\kappa}\mathbb{E}V(X(t),\mathcal{L}^1(X(t)),t)+\mathbb{E}{L}V(X_t,\mathcal{L}^1(X(t)),t)
\\&\leq -(\lambda-\bar{\kappa})\mathbb{E}V(X(t),\mathcal{L}^1(X(t)),t)<0.
\end{aligned}\end{align*}
Thanks to the continuity of the functions involved, one arrives at that there exists a sufficiently small $\bar h >0$ such that for any $0<h<\bar h$,
\begin{align}\label{eq5.15}
\begin{aligned}
&\bar{\kappa}\mathbb{E}V(X(s),\mathcal{L}^1(X(s)),s)+\mathbb{E}{L}V(X_s,\mathcal{L}^1(X(s)),s)
\leq0, \quad t\leq s\leq t+h.
\end{aligned}\end{align}
By Lemma \ref{le5.1} and \eqref{eq5.15}, we obtain that 
\begin{align}\label{eq5.14*}
\begin{aligned}
&e^{\bar{\kappa}(t+h)}\mathbb{E}V(X(t+h),\mathcal{L}^1(X(t+h)),t+h)-e^{\bar{\kappa}t}\mathbb{E}V(X(t),\mathcal{L}^1(X(t)),t)
\\&=\int_t^{t+h}e^{\bar{\kappa}s}[\bar{\kappa}\mathbb{E}V(X(s),\mathcal{L}^1(X(s)),s)+\mathbb{E}{L}V(X_s,\mathcal{L}^1(X(s)),s)]\mathrm{d}s\leq0.
\end{aligned}\end{align}
which implies that
\begin{align*}
e^{\bar{\kappa}(t+h)}\mathbb{E}V(X(t+h),\mathcal{L}^1(X(t+h)),t+h)\leq e^{\bar{\kappa}t}\mathbb{E}V(X(t),\mathcal{L}^1(X(t)),t).
\end{align*}
Therefore, $I(t+h)=I(t)$ for $0<h<\bar h$. Then one has $D^+I(t)=0$.
The proof is therefore complete.
\qed
 \begin{remark}\label{remark1}
If Theorem \ref{th3.2} holds for some $q\geq 2$, one can drive from Theorem \ref{T3} that for the initial value $\varsigma\in\mathbb{L}_{p}(C)(p\geq 4\vee l\vee q)$, the functional particle system \eqref{eq4.1} is exponentially stable in the $q$th  moment.
\end{remark}
\section{Example}\label{Sect.4*}
We provide an example to illustrate Theorem \ref{T3} and Theorem \ref{th3.2}.
\begin{expl}
\rm Consider the following scalar equation
\begin{align}\label{eq7.1}
\mathrm{d} y(t)&=\Big(-2 y(t)-3 y^3(t)+\frac{1}{4}\int_{-1/4}^0y(t+\theta)\mathrm{d}\theta+\frac{1}{2}\mathbb{E}^1y(t)\Big) \mathrm{d} t
\nonumber\\
&\quad+\frac{1}{2}( y(t)+\mathbb{E}^1y(t)) \mathrm{d} B_t
+\frac{1}{2}(y(t)+\mathbb{E}^1y(t)) \mathrm{d} B^0_t,
\end{align}
where $B_t$ and $B^0_t$ are two scalar Brownian motions defined on  $\left(\Omega^0, \mathcal{F}^0, \mathbb{P}^0\right)$ and $\left(\Omega^1, \mathcal{F}^1, \mathbb{P}^1\right)$, respectively. Here the initial data $y_0=\{y(\theta)=1, -1/4\leq \theta\leq 0\}$. It is easy to see that $y_0\in\mathbb{L}_4(C)$, and Assumption \ref{A1}, \ref{A2} and \ref{A3} hold. Set
$
V(x, \mu,t)=|x|^2+\int_{\mathbb{R}^d}|z |^2 \mu(\mathrm{d} z),
$
which satisfies \eqref{eq5.5} with $q=2$, $c_1=1$, $c_2=1$, $c_3=1$ and $c_4=1$.
By \cite[Vol-I, Example 3, p. 387]{CD2018}, one has
\begin{align*}
 &\partial_\mu\int_{\mathbb{R}^d}|z| ^2 \mu(\mathrm{d} z)(y)=2y;\quad\partial^2_\mu\int_{\mathbb{R}^d}|z| ^2 \mu(\mathrm{d} z)(y)=0;
\\&
 \partial_y\partial_\mu\int_{\mathbb{R}^d}|z| ^2 \mu(\mathrm{d} z)(y)=2;\quad\partial_x\partial_\mu\int_{\mathbb{R}^d}|z| ^2 \mu(\mathrm{d} z)(y)=0.
 \end{align*}
 By the definition of $LV$ in \eqref{lv}, we have
 \begin{align*}
{L}  V(\phi,\mathcal{L}^1(\phi(0)),t)
&=2\phi(0)\Big(-2 \phi(0)-3\phi^3(0)+\frac{1}{4}\int_{-1/4}^0\phi(\theta)\mathrm{d}\theta+\frac{1}{2}\mathbb{E}^1\phi(0)\Big)
\\
&\quad+2\Big(\frac{1}{2}(\phi(0)+\mathbb{E}^1\phi(0))\Big)^2
+2\mathbb{E}^1\Big[\phi(0)\Big(-2 \phi(0)-3\phi^3(0)\\
&\quad+\frac{1}{4}\int_{-1/4}^0\phi(\theta)\mathrm{d}\theta+\frac{1}{2}\mathbb{E}^1\phi(0)\Big)
+2\Big(\frac{1}{2}(\phi(0)+\mathbb{E}^1\phi(0))\Big)^2\Big].
\end{align*}
Using the Young inequality and the H\"older inequality, we compute that
\begin{align*}
{L}  V(\phi,\mathcal{L}^1(\phi(0)),t)
&\leq-4 |\phi(0)|^2-6|\phi(0)|^4+\frac{1}{4}|\phi(0)|^2+\frac{1}{4}\Big(\int_{-1/4}^0\phi(\theta)\mathrm{d}\theta\Big)^2
\\
&\quad+\frac{1}{2}|\phi(0)|^2+\frac{1}{2}|\mathbb{E}^1\phi(0)|^2+|\phi(0)|^2+|\mathbb{E}^1\phi(0)|^2
\\
&\quad+\mathbb{E}^1\Big[-4 |\phi(0)|^2-6|\phi(0)|^4+\frac{1}{4}|\phi(0)|^2+\frac{1}{4}\Big(\int_{-1/4}^0\phi(\theta)\mathrm{d}\theta\Big)^2
\\
&\quad+\frac{1}{2}|\phi(0)|^2+\frac{1}{2}|\mathbb{E}^1\phi(0)|^2+|\phi(0)|^2+|\mathbb{E}^1\phi(0)|^2\Big].
\\&\leq-\frac{9}{4}|\phi(0)|^2+\frac{1}{16}\int_{-1/4}^0|\phi(\theta)|^2\mathrm{d}\theta
\\&\quad+\frac{3}{4}\mathbb{E}^1|\phi(0)|^2+\frac{1}{16}\int_{-1/4}^0\mathbb{E}^1|\phi(\theta)|^2\mathrm{d}\theta.
\end{align*}
Taking expectation on both sides gives that
\begin{align}\label{eq6.2}
&\mathbb{E}{L}   V(\phi,\mathcal{L}^1(\phi(0)),t)
\leq
-\frac{3}{2}\mathbb{E}|\phi(0)|^2+\frac{1}{8}\int_{-1/4}^0\mathbb{E}|\phi(\theta)|^2\mathrm{d}\theta.
\end{align}
We can find a $\alpha>1$ such that ${3}/{4}-\alpha/64>0$. Therefore, for any $\phi\in\mathbb{L}_{t,2}(C)$ satisfying $\mathbb{E}(V(\phi(\theta), \mathcal{L}^1(\phi(\theta)),t+\theta))=2\mathbb{E}|\phi(\theta)|^2<\alpha\mathbb{E}(V(\phi(0), \mathcal{L}^1(\phi(0)),t))=2\alpha\mathbb{E}|\phi(0)|^2$ on $-\tau\leq \theta\leq0$, \eqref{eq6.2} yields
\begin{align*}
\mathbb{E}{L}   V(\phi,\mathcal{L}^1(\phi(0)),t)
&\leq -\Big(\frac{3}{2}-\frac{\alpha}{32}\Big)\mathbb{E}|\phi(0)|^2
\\&\leq
-\lambda\mathbb{E} V(\phi(0),\mathcal{L}^1(\phi(0)),t),
\end{align*}
where $$\lambda=\frac{3}{4}-\frac{\alpha}{64}>0.$$
\begin{figure}[htbp]
\centering
\includegraphics[height=5cm,width=15cm]{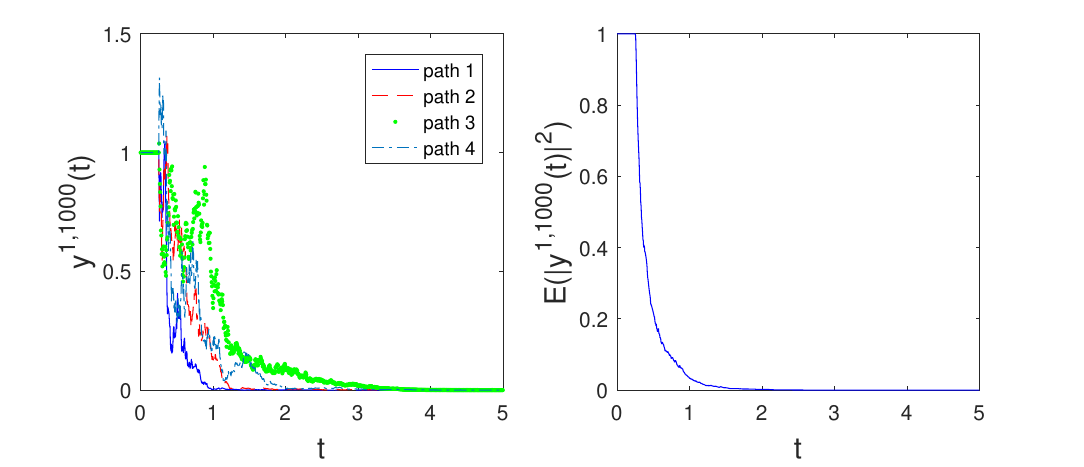}
\vspace{-1em}
\caption{Four sample paths of $y^{1,1000}(t)$ and the sample mean of $|y^{1,1000}(t)|^2$  of the equation \eqref{eq7.1*} for $t\in[0,5]$ with $100$ sample points and step size $\Delta=0.005$.}
\end{figure}
\!Subsequently, by Theorem \ref{th3.2}, the equation \eqref{eq7.1} is exponentially stable in mean square.
By virtue of Remark \ref{remark1}, the corresponding functional particle systems
\begin{align}\label{eq7.1*}
\mathrm{d} y^{k,N}(t)&=\Big(-2 y^{k,N}(t)-3 (y^{k,N}(t))^3+\frac{1}{4}\int_{-1/4}^0y^{k,N}(t+\theta)\mathrm{d}\theta+\frac{1}{2 N}\sum_{i=1}^Ny^{i,N}(t)\Big) \mathrm{d} t
\nonumber\\
&\quad+\frac{1}{2}\Big( y^{i,N}(t)+\frac{1}{N}\sum_{i=1}^Ny^{i,N}(t)\Big) \mathrm{d} B^k_t
+\frac{1}{2}\Big(y^{i,N}(t)+\frac{1}{N}\sum_{i=1}^Ny^{i,N}(t)\Big) \mathrm{d} B^0_t,
\end{align}
where $k=1,2,\cdots,N, t\geq0$, are exponentially stable in mean square.
Figure 1 depicts four sample paths of $y^{1,1000}(t)$ and the sample mean of $|y^{1,1000}(t)|^2$  of the equation \eqref{eq7.1*} with $N=1000$ for $t\in[0,5]$ with $100$ sample points and step size $\Delta=0.005$.
\end{expl}
\section{Conclusions}\label{Sect.5}
Considering past dependent phenomenon of the dynamic system, this paper focuses on the nonlinear MV-SFDEs with common noise \eqref{1.1}. The well-posedness is first investigated by a fixed-point argument. Secondly, the relationship between nonlinear MV-SFDEs with common noise and the corresponding mean-field functional particle systems is then established, encompassing the conditional propagation of chaos with an explicit bound on the convergence rate and the equivalence of stability.
Finally, the Razumikhin theorem is proven by utilizing the It\^o formula involved with state and measure that we developed for nonlinear MV-SFDEs with common noise. It provides an easily implementable criterion for the 
$q$th moment exponential stability of equation \eqref{1.1}.


\end{document}